\documentclass[letterpaper,10pt]{article}
\usepackage[margin=1in]{geometry}
\usepackage{graphicx}
\usepackage{import}
\usepackage[usenames,dvipsnames]{xcolor}
\usepackage[affil-it]{authblk}
\usepackage{float}
\usepackage[justification=centering]{caption}
\usepackage{algorithm,algpseudocode}
\usepackage{subcaption}
\usepackage{url}
\usepackage{amsmath}
\usepackage{amsthm}
\usepackage{amsfonts}
\usepackage{graphicx}
\usepackage[normalem]{ulem}
\usepackage[all]{xy}

\providecommand{\bR}{\mathbb{R}}
\providecommand{\bP}{\mathbb{P}}

\providecommand{\cF}{\mathcal{F}}
\providecommand{\cG}{\mathcal{G}}

\DeclareMathOperator{\spn}{span}

\DeclareMathOperator{\Gr}{Gr}
\DeclareMathOperator{\SO}{SO}

\theoremstyle{plain}
\newtheorem*{mthm}{Main Theorem}
\newtheorem{proposition}{Proposition}
\newtheorem{lemma}{Lemma}
\newtheorem{theorem}{Theorem}
\newtheorem{corollary}{Corollary}
\newtheorem{problem}{Problem}

\theoremstyle{definition}
\newtheorem{definition}{Definition}
\newtheorem{example}{Example}
\newtheorem{remark}{Remark}

\numberwithin{equation}{section}

\numberwithin{proposition}{section}
\numberwithin{lemma}{section}
\numberwithin{theorem}{section}
\numberwithin{corollary}{section}

\numberwithin{definition}{section}

\numberwithin{example}{section}
\numberwithin{remark}{section}

\DeclareMathOperator{\dO}{O}

\newcommand{\E}{(\bP^3)^{\binom{N}{2}}}

\DeclareMathOperator{\sign}{sign}
\begin{document}
\title{\vspace{-6ex}Semi-algebraic geometry of common lines}
\author{\vspace{-1.5ex}David Dynerman\thanks{Department of Mathematics, University of Wisconsin -- Madison \texttt{dynerman@math.wisc.edu}}}
\date{}
\maketitle
\begin{abstract}
Cryo-electron microscopy is a technique in structural biology for discovering/determining the 3D structure of small molecules. A key step in this process is detecting common lines of intersection between unknown embedded image planes. We intrinsically characterize such common lines in terms of the unembedded geometric data detected in experiments. We show these common lines form a semi-algebraic set, i.e., they are defined by polynomial equalities and inequalities. These polynomials are low degree and, using techniques from spherical geometry, we explicitly derive them in this paper.
\end{abstract}

\section{Introduction}
\label{sect-intro}
Cryo-electron microscopy (cryo-EM) is a technique used to discover the structure of small molecules, usually proteins in the context of structural biology research~\cite{wang}. 
\begin{figure}[H]
\begin{center}
\includegraphics[width=4in]{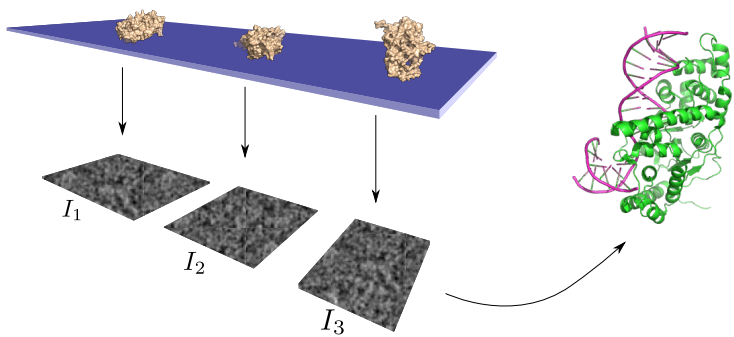}
\vspace{-14ex}
\end{center}
\caption{Cryo-EM reconstructs a 3D structure from noisy 2D images $I_1, \ldots, I_N$.}
\label{fig-overview}
\end{figure}
A basic outline of cryo-EM is presented in [Figure~\ref{fig-overview}]. First, a sample is prepared by freezing many different copies of the molecule in a thin layer of ice. The sample is then illuminated by a stream of electrons which are detected by cameras and produce $N$ noisy 2D cryo-EM images $I_1, \ldots, I_N$. The primary goal is to reconstruct the 3D structure of the molecule from the 2D images that are acquired. For a more detailed overview, see~\cite[Section 1]{singer-cryo-num}.
\begin{problem}[Reconstruction Problem: Structural Biology]
\label{prob-bio}
Given $N$ two dimensional experimental cryo-EM images $I_1, \ldots, I_N$, reconstruct a three-dimensional model of the original molecule.
\end{problem}

\subsection{Mathematical Model}
\begin{figure}
\centering
\begin{subfigure}[b]{0.4\textwidth}
\centering
\includegraphics[width=1.5in]{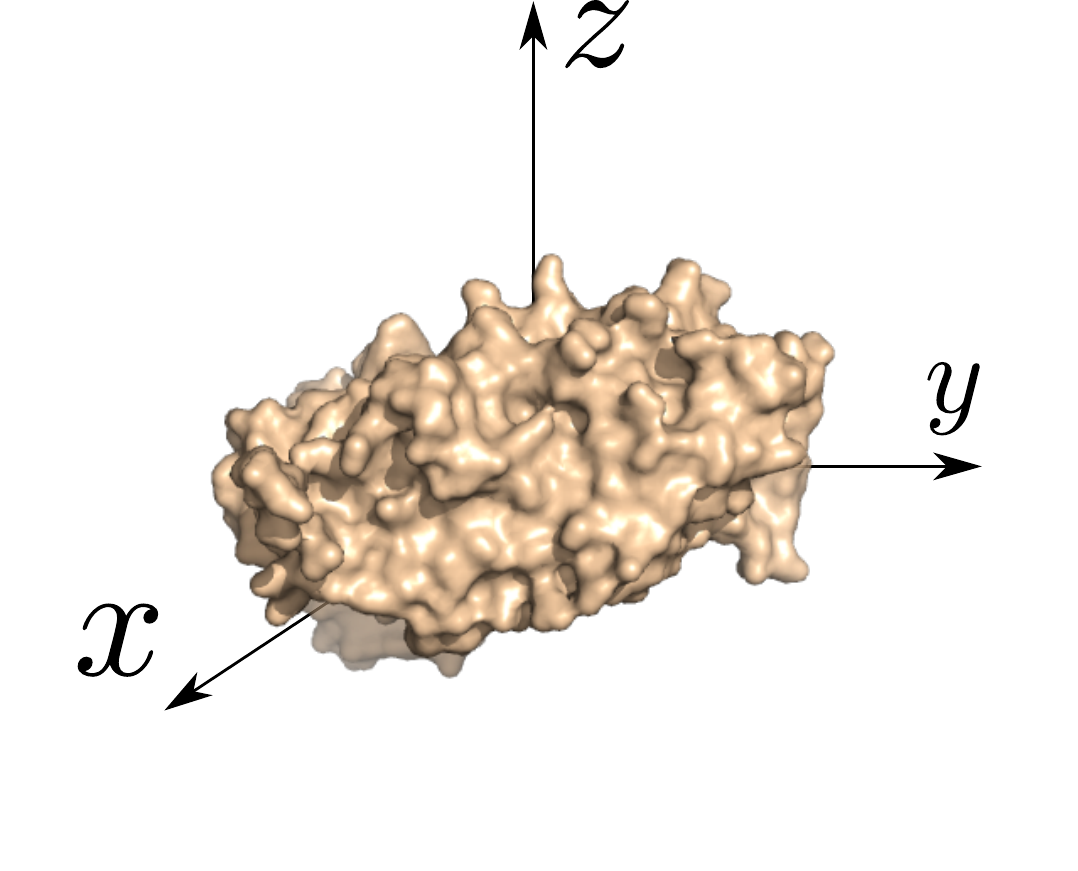}
\caption{Molecule $\phi: \bR^3 \to \bR$.}
\label{fig-model-phi}
\end{subfigure}
\quad
\begin{subfigure}[b]{0.4\textwidth}
\centering
\includegraphics[width=1.5in]{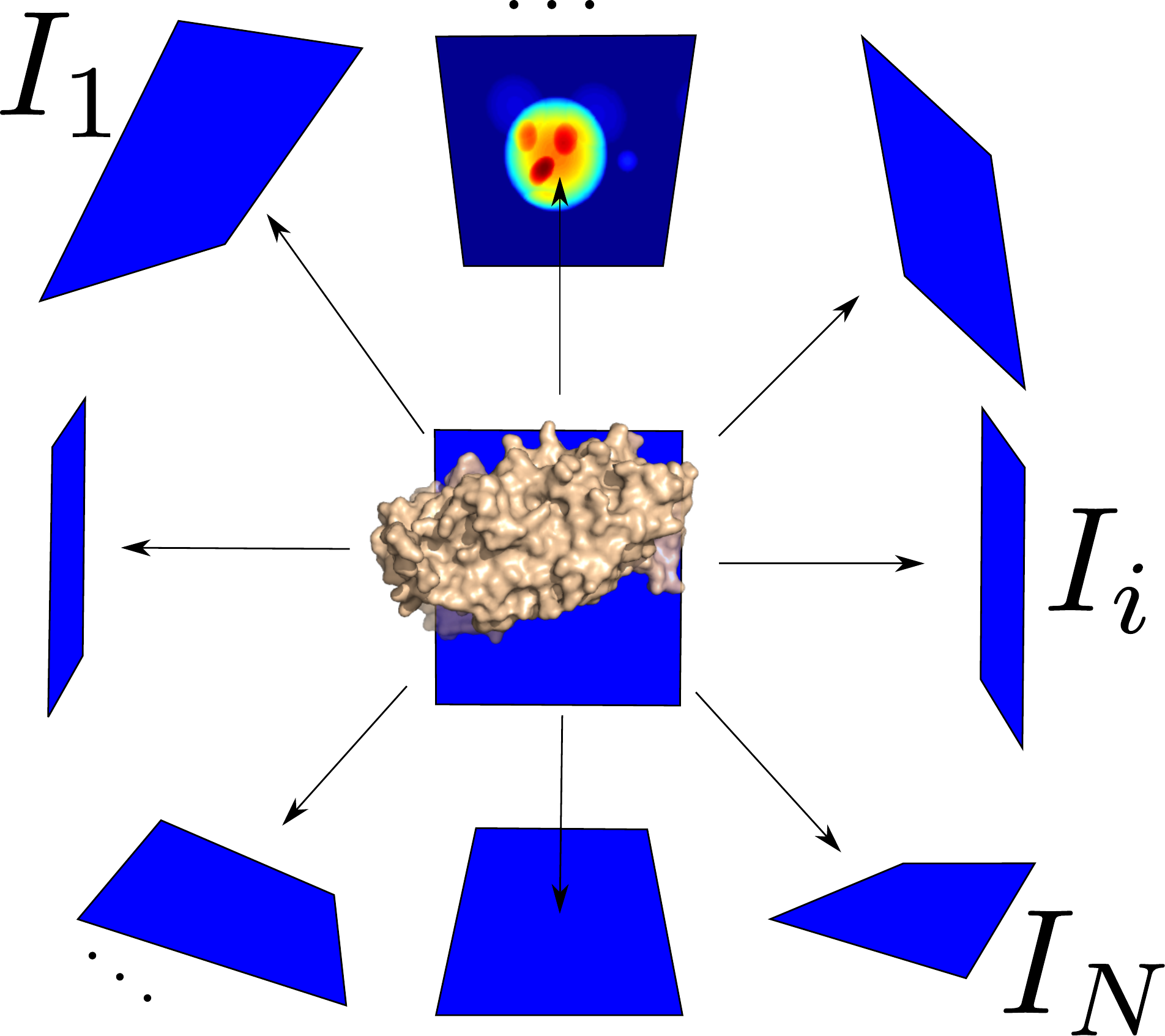}
\caption{Many images of a single molecule.}
\label{fig-model-het}
\end{subfigure}

\vspace{0.25in}

\begin{subfigure}[b]{0.4\textwidth}
\centering
\includegraphics[width=1.5in]{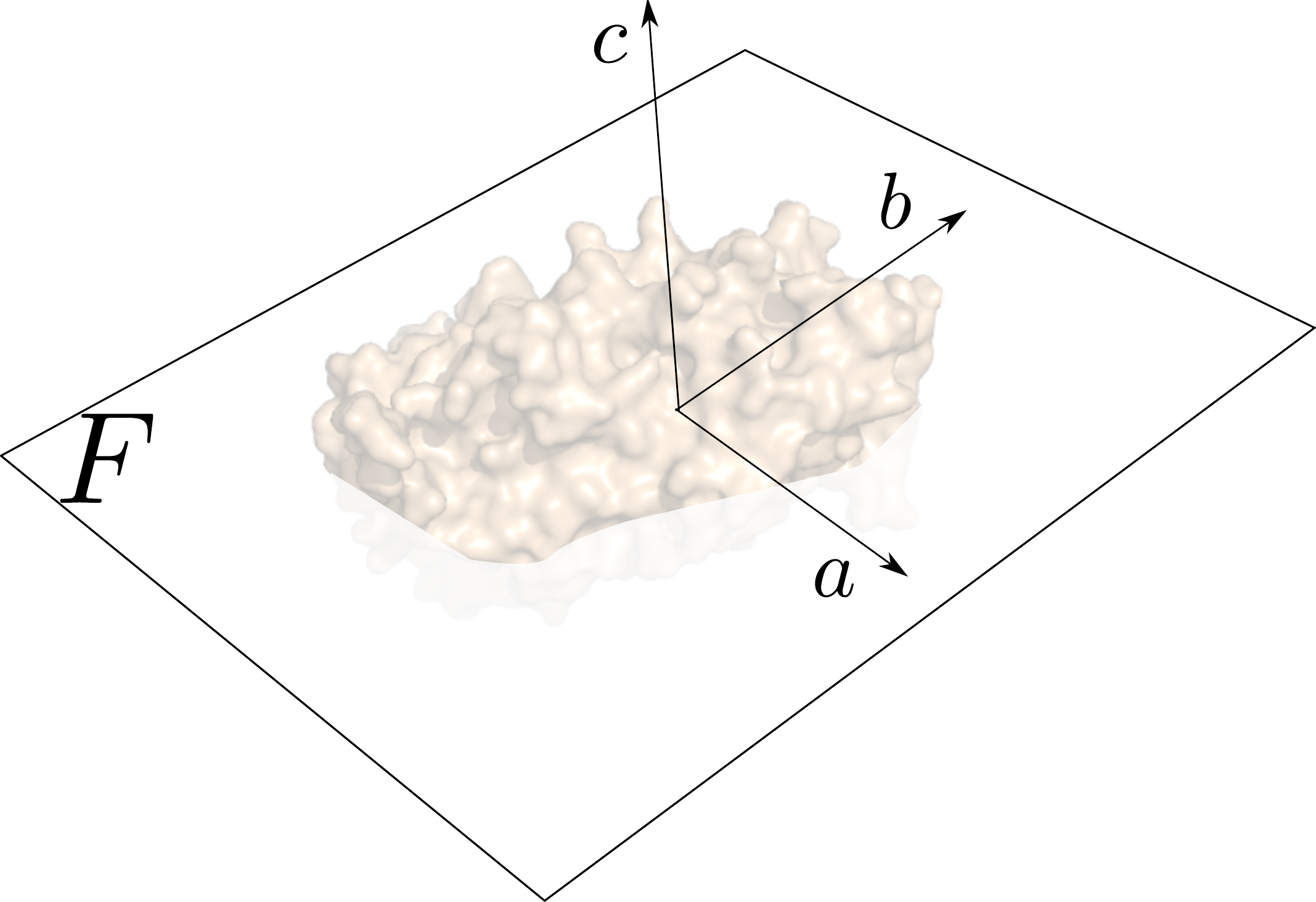}
\caption{Microscope orientation $F$.}
\label{fig-model-frame}
\end{subfigure}
\quad
\begin{subfigure}[b]{0.4\textwidth}
\centering
\includegraphics[width=1.5in]{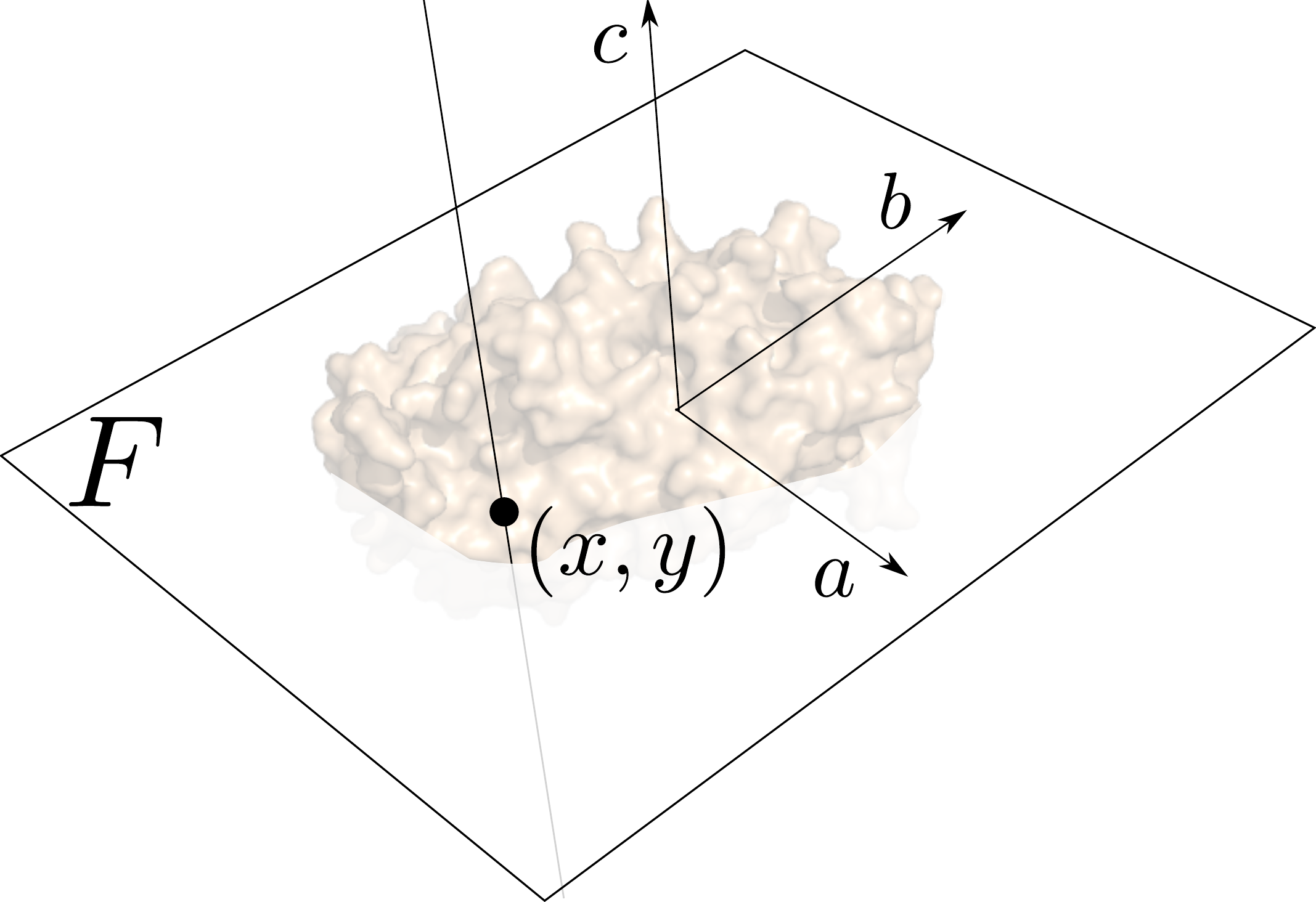}
\caption{Cryo-EM image from orientation $F$.}
\label{fig-model-image}
\end{subfigure}
\caption{Cryo-EM mathematical model.}
\end{figure}

We briefly describe the mathematical model for cryo-EM, following~\cite[Section 0]{singer-cryo}. We work in the three dimensional space $\bR^3$ equipped with the usual inner product. The molecule is modeled by a function $\phi: \bR^3 \to \bR$ that represents its electronic density at various spatial locations [Figure~\ref{fig-model-phi}]. An actual cryo-EM experiment obtains a single image of many copies of the molecule, but we instead assume that each image is a picture of the same molecule from different microscope orientations [Figure~\ref{fig-model-het}]. To model a microscope orientation we use the following concept.
\begin{definition}
A \uline{frame} $F$ for $\bR^3$ is an ordered orthonormal basis $(a,b,c)$ such that the determinant of the matrix $[a\, b\, c]$ is $+1$, or, equivalently, that $c = a \times b$, where $\times$ is the standard cross product on $\bR^3$. 
\end{definition}
\begin{remark}
A frame $F$ for $\bR^3$ is uniquely determined by the vectors $(a,b)$. For the rest of the paper we identify frames $(a,b,c)$ with pairs of orthonormal vectors $(a,b)$. 
\end{remark}
For us a microscope orientation is a frame $F = (a,b)$. We think of the span of the vectors $a$ and $b$ as the embedded image plane of this orientation, and the vector $c = a \times b$ as the ``viewing'' direction [Figure~\ref{fig-model-frame}].
 
A cryo-EM experiment produces $N$ images which we denote $I_1, \ldots, I_N$ --- see [Figure~\ref{fig-model-het}]. We will write $F_i = (a_i, b_i)$ for the microscope orientation of image $I_i$. The embedded image plane spanned by $a_i, b_i$ can be canonically identified with the plane $P_i = \bR^2$. We think of $P_i$ as the \emph{unembedded} image plane of $I_i$. We model the image $I_i$ as a real valued function on $P_i = \bR^2$. The value of the image $I_i$ at the point $(x,y)$ is the integral of $\phi$ along a line perpendicular to the embedded image plane $\spn \{ a_i, b_i\}$ --- see [Figure~\ref{fig-model-image}] and [Equation~\ref{eqn-image}]. This is the \emph{X-ray transform} of $\phi$ onto the frame $F_i$, given by
\begin{equation}
\label{eqn-image}
\begin{aligned}
I_i&: P_i = \bR^2 \to \bR, \\
I_i(x,y) &= \int_{-\infty}^\infty \phi(xa_i + yb_i + zc_i) dz,
\end{aligned}
\end{equation}
where $c_i = a_i \times b_i$. As in~\cite{singer-cryo}, to solve this reconstruction problem we assume that the X-ray projections $I_i$ and $I_j$ of $\phi$ from different microscope orientations $F_i$ and $F_j$ are different. This is equivalent to requiring the molecule $\phi$ to admit no non-trivial symmetry as a function on $\bR^3$. 

In terms of this mathematical model, the goal of cryo-EM reconstruction [Problem~\ref{prob-bio}] becomes to recover the function $\phi$ from the $N$ X-ray projections $I_1, \ldots, I_N$. A commonly used approach for this problem is to first recover the $N$ projection orientations $F_1, \ldots, F_N$~\cite[Section 0.1]{singer-cryo}. Note that the detected image $I_i$ is a function on the plane $P_i = \bR^2$, and a cryo-EM experiment does not directly provide information about the microscope orientation $F_i$ used to compute $I_i$. 

Once the original microscope orientations are known, the unembedded image data $I_1, \ldots, I_N$ can be placed in the original positions from where these X-ray projections were computed. Then the X-ray transform can be inverted to yield an approximation of $\phi$. Thus, although the ultimate goal is to solve [Problem~\ref{prob-bio}], we instead discuss solutions to the following problem.
\begin{problem}[Reconstruction Problem: Microscope Orientations]
\label{prob-frames}
Given $N$ X-ray projections $I_1, \ldots, I_N$ of a molecule $\phi: \bR^3 \to \bR$, computed from the $N$ unknown microscope orientations $F_1, \ldots, F_N$, recover these orientations up to global rotation.
\end{problem}
By ``up to global rotation'' we mean that instead of recovering $\phi$ exactly, we might recover a rotated version of $\phi$ by a transformation $R$ in the group $\dO(3)$ of all $3 \times 3$ rotation matrices. Rotational ambiguity in the reconstructed molecule is not a problem, since recovering a rotated version of the molecule is as good as recovering the original. It may be the case that $R$ is an improper rotation, i.e., $\det R = -1$, in which case the recovered version will have the opposite chirality of the original molecule which is not desirable. However, other techniques exist to resolve this chiral ambiguity, so this also does not pose a problem. 

\subsection{Common Lines and Reconstruction}
\label{sect-common}
\begin{figure}[ht]
\begin{center}
\includegraphics[width=4in]{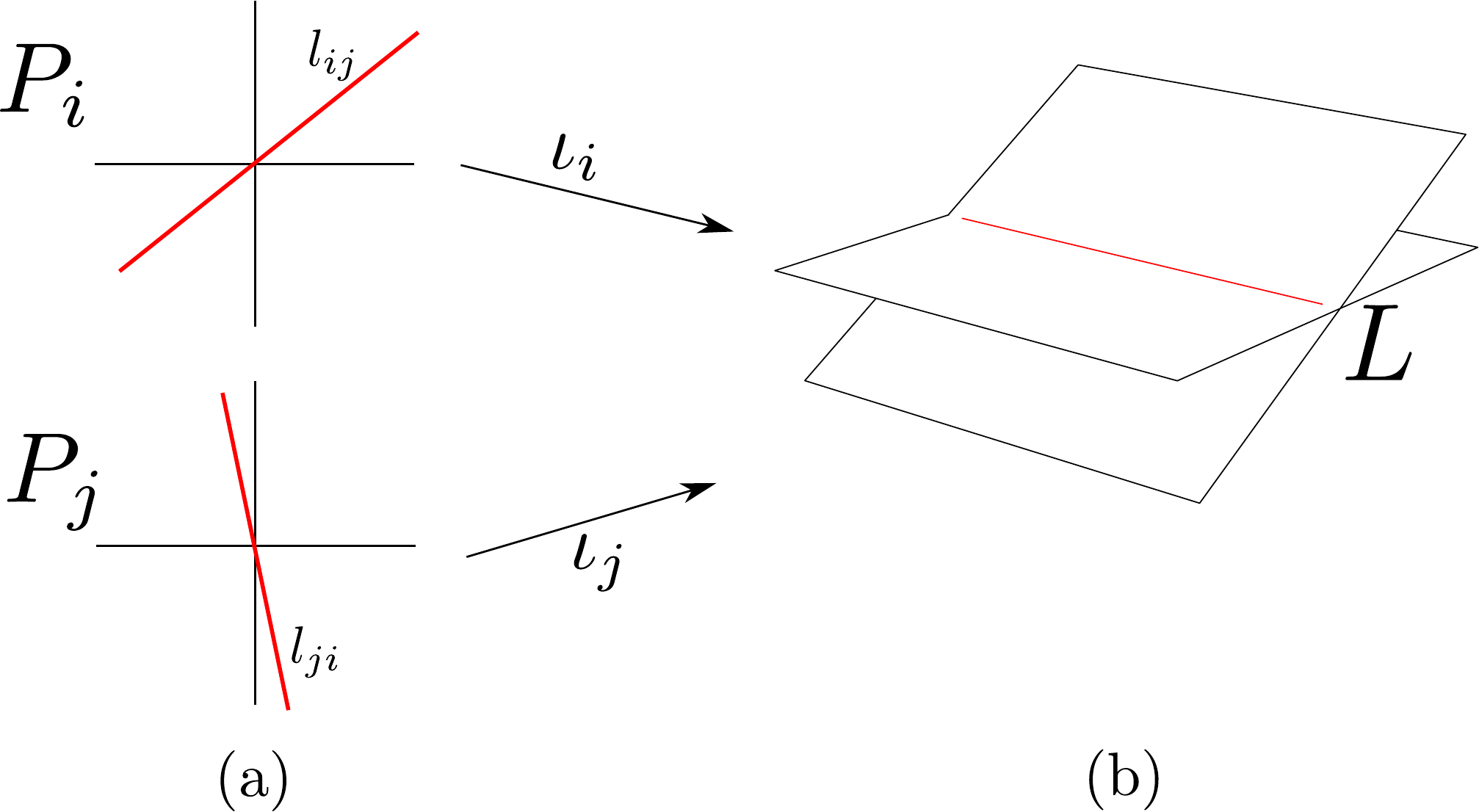}
\end{center}
\caption{Common line of $F_i$ and $F_j$.}
\label{fig-common-ij}
\end{figure}

One approach for solving [Problem~\ref{prob-frames}] is to exploit common lines of intersection between the embedded image planes, which we now describe. A cryo-EM experiment produces images $I_i$ and $I_j$ from orientations $F_i = (a_i, b_i)$ and $F_j = (a_j, b_j)$. These frames define isometric embeddings $\iota_i$ and $\iota_j$ [Figure~\ref{fig-common-ij}] of the unembedded image planes $P_i$ and $P_j$ into $\bR^3$, given by
\[
\iota_i(x,y) = xa_i + yb_i, \quad \quad \iota_j(x,y) = xa_j + yb_j.
\]
The images are functions on $P_i$ and $P_j$, and we know that they were obtained as X-ray projections onto the unknown embedded image planes $\iota_i(P_i)$ and $\iota_j(P_j)$ [Figure~\ref{fig-common-ij}b]. As in~\cite{vanheel,singer-cryo-num} we assume that the unknown microscope orientations are sampled uniformly from the space of all frames. This implies that the planes $\iota_i(P_i)$ are distinct, and, further, that each pair of these planes intersects in a different line. Such a configuration of frames is called \emph{generic}.

The embedded image planes $\iota_i(P_i)$ and $\iota_j(P_j)$ intersect in a line $L$ [Figure~\ref{fig-common-ij}b], and this line corresponds to the unembedded lines $\ell_{ij} \subset P_i$ and $\ell_{ji} \subset P_j$. Since these unembedded lines both came from $L \subset \bR^3$ we have a natural choice of isometry\footnote{Note that there are only two possible isometries between $\ell_{ij}$ and $\ell_{ji}$.} $\psi_{ij}: \ell_{ij} \to \ell_{ji}$. Proceeding in this fashion, the $N$ microscope orientations $F_1, \ldots, F_N$ produce\footnote{For positive integers $N > k$, ``$N$ choose $k$'' is the integer $\binom{N}{k} = \frac{N!}{k!(N-k)!}$. This is the number of ways to choose $k$ distinct numbers from $\{ 1, \ldots, N \}$.} $\binom{N}{2}$ common line pairs $\{ (\ell_{ij}, \ell_{ji}, \psi_{ij}) \}$. This is the \emph{common lines data realized} by the frames $F_1, \ldots, F_N$. It will be useful for us to distinguish such common lines data obtained from frames.
\begin{definition}
A \uline{common line pair} for $P_i$ and $P_j$ is a pair of lines $\ell_{ij} \subset P_i$ and $\ell_{ji} \subset P_j$, together with a choice of isometry $\psi_{ij}: \ell_{ij} \to \ell_{ji}$. A collection of common line pairs $\{ (\ell_{ij}, \ell_{ji}, \psi_{ij}) \}$, for every $P_i$ and $P_j$, is \uline{common lines data} for $P_1, \ldots, P_N$. We say common lines data is \uline{valid} if it is realized by some generic frames $F_1, \ldots, F_N$.
\end{definition}

Despite the fact that common lines data is information in the unembedded planes $P_i$, it is a fact that valid common lines data determines its realizing frames, up to global rotation. Further, algorithms have long been known (e.g.~\cite[Section 2.1]{vanheel}) that recover a set of realizing frames from valid common lines data. 

This is relevant to cryo-EM reconstruction, because although the microscope orientations are unknown, it is possible to detect the common lines data the orientations realize from the images $I_1, \ldots, I_N$~\cite[Equation 2.3]{singer-cryo-num}. Thus we have the following \emph{common lines approach} for the cryo-EM reconstruction problem [Problem~\ref{prob-frames}]. We first detect the common lines data realized by the unknown microscope orientations. Next, from the valid common lines data we reconstruct a set of realizing frames. Since valid common lines data determines its realizing frames up to global rotation, the reconstructed frames are related to the original microscope orientations by a global rotation, and so in principle one has solved the reconstruction problem. 

\subsection{Angular Reconstruction}
\label{sect-vanheel}
In this section we describe the \emph{angular reconstruction algorithm}, due to van Heel~\cite{vanheel}, and also independently Vainshtein and Goncharov~\cite{vainshtein}, which recovers a set of realizing frames from valid common lines data. 

\begin{figure}[ht]
\begin{center}
\begin{minipage}[c]{0.3\textwidth}
\centering\includegraphics[width=1.25in]{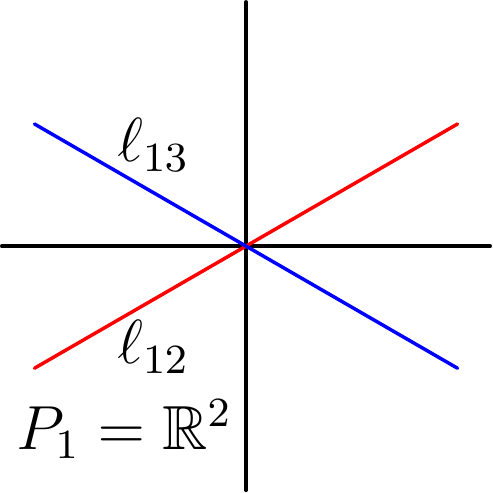}
\end{minipage}
\begin{minipage}[c]{0.3\textwidth}
\centering\includegraphics[width=1.25in]{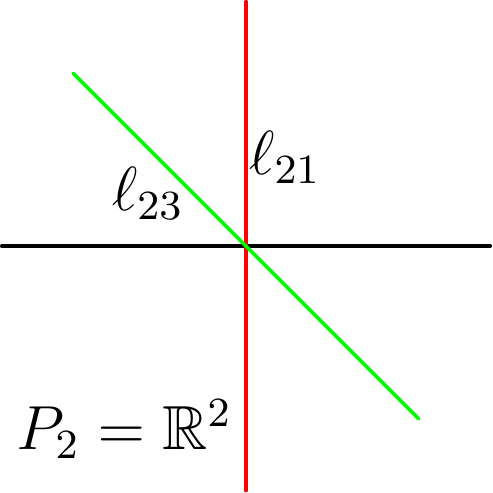}
\end{minipage}
\begin{minipage}[c]{0.3\textwidth}
\centering\includegraphics[width=1.25in]{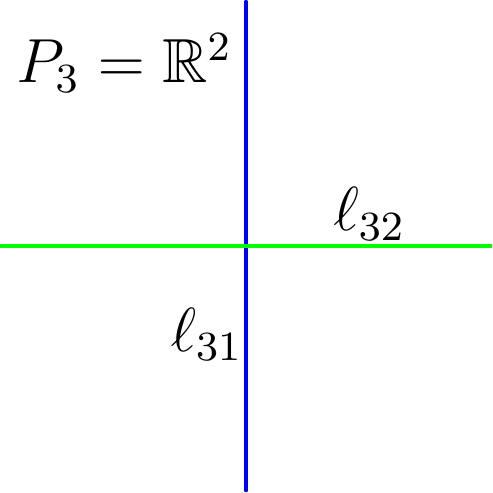}
\end{minipage}
\end{center}
\caption{Common lines data for $N = 3$ planes.}
\label{fig-vanheel-input}
\end{figure}

\begin{figure}[ht]
\begin{center}
\begin{subfigure}[t]{0.225\textwidth}
\centering\includegraphics[width=1.5in]{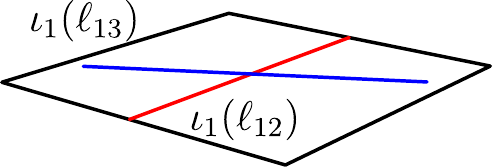}
\caption{Place $P_1$.}
\label{fig-vanheel-1}
\end{subfigure}
\begin{subfigure}[t]{0.225\textwidth}
\centering\includegraphics[width=1.5in]{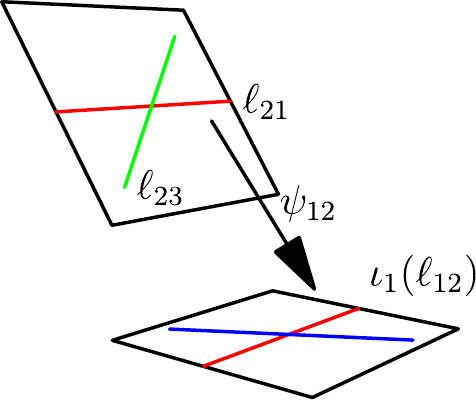}
\caption{Dock $P_2$ via $\psi_{12}$.}
\label{fig-vanheel-2}
\end{subfigure}
\begin{subfigure}[t]{0.225\textwidth}
\centering\includegraphics[width=1.5in]{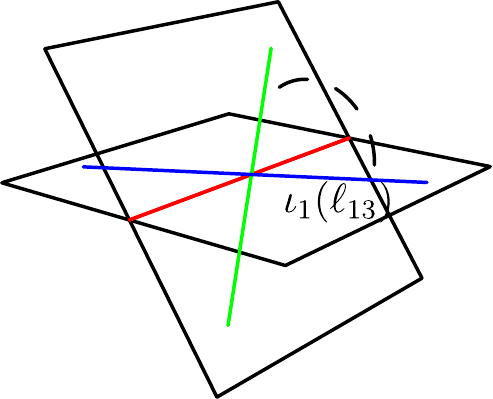}
\caption{Rotational ambiguity for $\iota_2(P_2)$.}
\label{fig-vanheel-3}
\end{subfigure}
\begin{subfigure}[t]{0.225\textwidth}
\centering\includegraphics[width=2in]{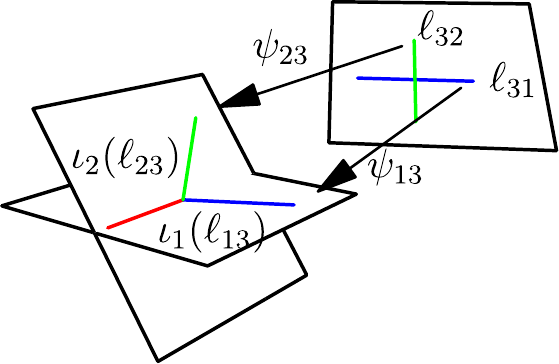}
\caption{Resolve by docking $P_3$ to $\iota_1(\ell_{13})$.}
\label{fig-vanheel-4}
\end{subfigure}
\end{center}
\caption{Angular reconstruction.}
\label{fig-vanheel}
\end{figure}

Our input is valid common lines data $\{ (\ell_{ij}, \ell_{ji}, \psi_{ij}) \}$ for $P_1, \ldots, P_N$ [Figure~\ref{fig-vanheel-input}]. Note that recovering a frame $F_i$ is equivalent to recovering the embedding $\iota_i$ of $P_i$, which will be easier to visualize. Since we are only reconstructing up to global rotation, the first step is to embed $P_1$ in an arbitrary position in $\bR^3$ [Figure~\ref{fig-vanheel-1}]. Next, we use the isometry $\psi_{12}$ between $\ell_{12}$ and $\ell_{21}$ to dock $P_2$ to $\iota_1(P_1)$ [Figure~\ref{fig-vanheel-2}]. This docking is ambiguous [Figure~\ref{fig-vanheel-3}] since we are free to rotate $\iota_2(P_2)$ about its line of intersection with $\iota_1(P_1)$. We resolve this ambiguity by docking $P_3$ with $\iota_1(P_1)$ and matching up $\ell_{23}$ and $\ell_{32}$ in $\iota_2(P_2)$ and $\iota_3(P_3)$ [Figure~\ref{fig-vanheel-4}]. We continue in this fashion, docking each subsequent plane $P_i$ with $\iota_1(P_1)$ and resolving the rotational ambiguity by comparing against the remaining frames. 

\subsection{Noise and valid common lines data}
\label{sect-noise}
\begin{figure}[h]
\centering
\includegraphics[width=1in]{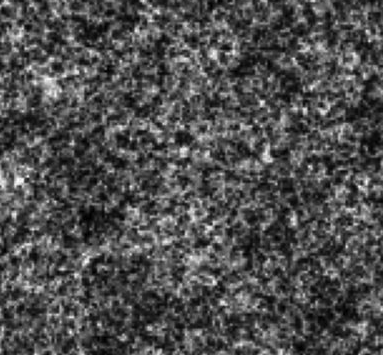}
\caption{Raw cryo-EM image~\cite{henderson}.}
\label{fig-noisy}
\end{figure}

We discussed in [Section~\ref{sect-common}] that valid common lines data determines its realizing frames up to global rotation. Common lines based approaches for cryo-EM reconstruction [Problem~\ref{prob-frames}] assume that we can accurately detect the valid common lines realized by the unknown microscope orientations. Unfortunately cryo-EM images are very noisy [Figure~\ref{fig-noisy}], so we cannot expect to correctly identify common lines data.

Misdetected common lines pose a problem because they lead to inconsistencies when attempting to recover realizing frames. For example, in [Figure~\ref{fig-vanheel}] we resolved the ambiguity of $\iota_2(P_2)$ by docking $P_3$ to $\iota_1(P_1)$ and using the common lines $l_{23}$ and $l_{32}$ [Figure~\ref{fig-vanheel-3}]. However, we could have equally well resolved the ambiguity of $\iota_2(P_2)$ by docking $P_4$ and using the common lines $l_{24}$ and $l_{42}$. Thus, if we, for example, incorrectly identify the common lines in $P_4$ we will have two contradictory embeddings $\iota_2(P_2)$ with no obvious way of determining which is correct. 

More generally, the angular reconstruction algorithm makes many choices: for example which plane to begin reconstruction with, and how to resolve docking ambiguities. The final reconstructed frames depend on all these choices. By definition valid common lines data is precisely the data which has a single consistent (up to global rotation) set of realizing frames. The development of common lines reconstruction algorithms that are robust to this kind of error is an active area of research. 

\subsection{Our Results}
\label{sect-results}
We wish to understand the set $C_N$ of all valid common lines data for $N$ planes $P_1, \ldots, P_N$. First, we derive necessary and sufficient conditions for common lines data to be valid. These conditions are polynomial equations and inequalities, which means that $C_N$ is a \emph{semi-algebraic set}, and allows us to study $C_N$ as a geometric space. In particular, we compute the dimension of $C_N$, and show that there is a geometric bijection between $C_N$ and the space of generic frames, up to global rotation. 
\begin{mthm}
The set $C_N$ of all valid common lines data for $N$ frames is a $3N-3$ dimensional semi-algebraic subset of the $2\binom{N}{2}$ dimensional space of all common lines data, and is homeomorphic to the space of $N$ generic frames modulo $\dO(3)$. The defining equations for $C_N$ are given by $\binom{N}{3}$ polynomial inequalities arising from the spherical triangle inequalities and $6\binom{N}{4}$ polynomial equalities arising from the spherical law of cosines.
\end{mthm}
The meaning of this theorem is as follows. As we discussed in [Section~\ref{sect-common}], one way to obtain valid common lines data is from the embedded frames $F_1, \ldots, F_N$. The theorem provides an intrinsic definition of this valid common lines data, namely, the defining polynomials for $C_N$. This is a definition for valid common lines only in terms of the data $\{ (l_{ij}, l_{ji}, \psi_{ij}) \}$ on unembedded planes $P_1, \ldots, P_N$, and without reference to any embedded frames $F_1, \ldots, F_N$.

\begin{figure}[h]
\begin{center}
\begin{subfigure}[t]{0.3\textwidth}
\centering\includegraphics[width=1.5in]{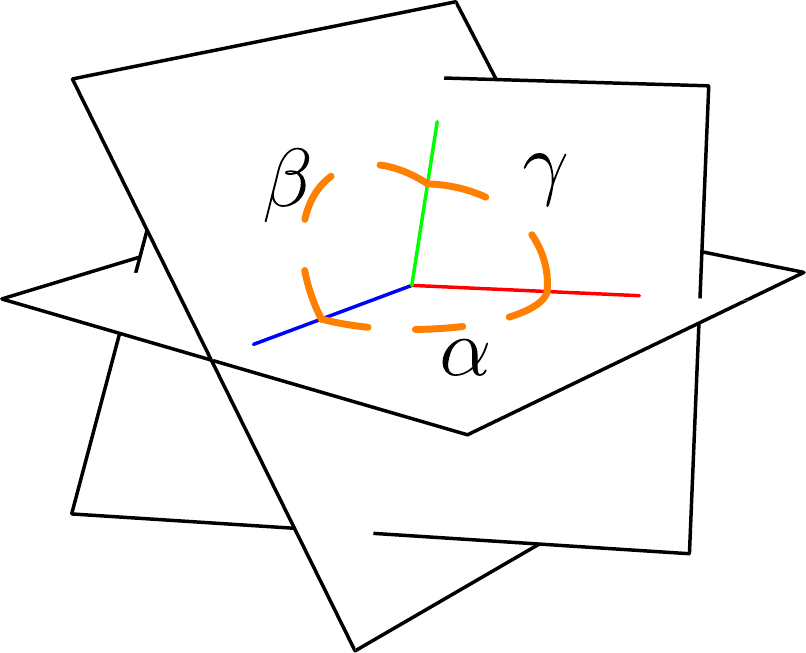}
\end{subfigure}
\begin{subfigure}[t]{0.3\textwidth}
\centering\includegraphics[width=1.25in]{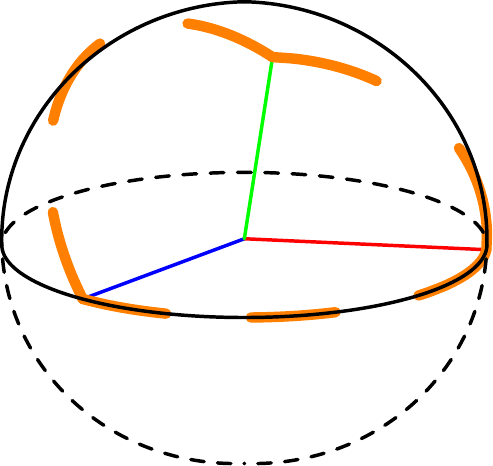}
\end{subfigure}
\end{center}
\caption{Common lines in $P_1$, $P_2$ and $P_3$ determine a spherical triangle.}
\label{fig-vanheel-tri}
\end{figure}

We briefly describe the idea behind our proofs. Suppose we have valid common lines data 
\begin{equation}
\label{eqn-common-3}
\{ (\ell_{12}, \ell_{21}, \psi_{12}), (\ell_{13}, \ell_{31}, \psi_{13}), (\ell_{23}, \ell_{32}, \psi_{23}) \}.
\end{equation}
The angles between these unembedded common lines determine a spherical triangle\footnote{By a spherical triangle we mean the data of $3$ points on the unit sphere in $\bR^3$, together with geodesic arcs on the sphere joining these vertices} [Figure~\ref{fig-vanheel-tri}], and so the angles $\alpha$ between $\ell_{12}$ and $\ell_{13}$, $\beta$ between $\ell_{21}$ and $\ell_{23}$ and $\gamma$ between $\ell_{31}$ and $\ell_{32}$ must satisfy the \emph{spherical triangle inequalities}. These inequalities are analogs of the plane triangle inequality, i.e. necessary and sufficient conditions for a spherical triangle to exist with the specified edge lengths. In other words, a necessary condition for common lines data to be valid is that it satisfy such triangle inequalities. In fact, we will see that for $N = 3$ having the common lines [Equation~\ref{eqn-common-3}] satisfy the spherical triangle inequalities is sufficient for the data to be valid [Proposition~\ref{prop-three}]. 
We prove our results for $N > 3$ by similarly appealing to spherical trigonometry. Specifically, given common lines data $\{ \ell_{ij}, \ell_{ji}, \psi_{ij} \}$ for $N$ planes, we require that for each triple $1 \leq i < j < k \leq N$ the common lines data $(\ell_{ij}, \ell_{ji}, \psi_{ij})$, $(\ell_{ik}, \ell_{ki}, \psi_{ik})$ and $(\ell_{jk}, \ell_{kj}, \psi_{jk})$ satisfy the spherical triangle inequalities. Now, reducing to the $N = 3$ case gives us realizing embeddings $\iota_i$, $\iota_j$, $\iota_k$ for \emph{each} triple $(i,j,k)$ of indices. To reconstruct a collection of $N$ consistent frames, all these triple reconstructions must be compatible. We show that this compatibility condition is a polynomial condition arising from the spherical law of cosines. These defining equations are given by polynomials which are explicitly derived and listed in [Section~\ref{sect-eqns}]. 

\subsection{Future Work}
Thinking of valid common lines data in geometric terms provides some insight about inconsistencies during reconstruction due to noise. The space of all common lines data has dimension $N(N-1)$, and, since valid common lines are homeomorphic to the space of $N$ frames up to global rotation, we have that the dimension of $C_N$ is $3N - 3$. Since $C_N$ is a space of small dimension in the ambient space, it follows that the reconstruction inconsistencies described in [Section~\ref{sect-noise}] are guaranteed to occur. In effect the most basic version of the angular reconstruction algorithm reconstructs the microscope orientations $F_1, \ldots, F_N$ using only $2N - 3$ out of the $\binom{N}{2}$ common line pairs, and arbitrarily ignores inconsistencies within these pairs. The set $C_N$ is precisely the set of common lines data for which this algorithm will produce the same output regardless of which common line pairs are used, but as described above we do not expect experimental data to lie in $C_N$.

Developing common lines reconstruction algorithms that are robust to noise is an active area of research. We are interested in exploring a geometric approach to noise reduction, which we briefly describe. In principle noisy experimental data $\{ (l'_{ij}, l'_{ji}, \psi'_{ij}) \}$ that lies outside of $C_N$ ``came from'' some noiseless valid common lines data in $C_N$. Since the set $C_N$ is the set of solutions of a system of polynomials, it is theoretically possible to project noisy common lines to the set of noiseless common lines $C_N$ via constrained polynomial optimization. We are interested in developing effective projection algorithms along these lines to reduce the impact of noise in reconstruction.

\section{Defining Equations}
\label{sect-eqns}
We proceed to describe in detail the results in [Main Theorem]. We derive necessary and sufficient conditions for common lines data to be valid. These will be explicit polynomial equations and inequalities only in the unembedded information $\{ (l_{ij}, l_{ji}, \psi_{ij}) \}$, and will provide an intrinsic definition for valid common lines without reference to the frames $F_1, \ldots, F_N$. We only consider the case when $N \geq 3$.

\subsection{Projective coordinates}
To obtain coordinates for $C_N$ it will be convenient for us to work with \emph{projective coordinates}, which we briefly review. Suppose $V$ is a vector space and $\ell$ is a line in $V$ through the origin. We can represent $\ell$ by choosing any non-zero vector $v \in \ell$. In other words, lines can be identified with equivalence classes of vectors. We denote the equivalence class of a vector $v$ by $[v]$, and by definition $[v] = [w]$ if and only if the vectors $v = \lambda w$, for $\lambda \neq 0$. The space of all lines through the origin in $V$ is the projective space $\bP(V)$. If $V = U \times W$ and $(u,w) \in V$, then we write $[u:w]$ for the corresponding class in $\bP(U \times W)$.

\subsection{Coordinates for Common Lines}
\label{sect-coords}
Suppose now that $(\ell_{ij}, \ell_{ji}, \psi_{ij})$ is a common line pair for $P_i$ and $P_j$. Choose a vector $v_{ij} = (x_{ij}, y_{ij})$ on the line $\ell_{ij} \subset P_i$, and consider the pair $(v_{ij}, \psi_{ij}(v_{ij})) \in P_i \times P_j$. Note that different choices of a vector along $\ell_{ij}$ will simply scale $(v_{ij}, \psi_{ij}(v_{ij}))$ by a non-zero multiple, so the projective pair $[v_{ij} : \psi_{ij}(v_{ij})]$ in $\bP(P_i \times P_j)$ is uniquely determined. 

Conversely, if $[v_{ij} : v_{ji}]$ satisfies $\| v_{ij} \|^2 = \| v_{ji} \|^2$, we obtain a common line pair $(\spn \{ v_{ij} \}, \spn \{ v_{ji} \}, \psi_{ij})$, where $\psi_{ij}$ is the unique isometry that sends $v_{ij} \mapsto v_{ji}$. Note that we obtain the same common line pair regardless of which representing vectors $(v_{ij}, v_{ji})$ we choose.

Thus, from now on we identify common line pairs with elements $[v_{ij} : v_{ji}] \in \bP(P_i \times P_j)$ satisfying $\| v_{ij} \|^2 = \| v_{ji} \|^2$. We also apply this identification to common lines data:
\begin{remark}
We identify common lines data for $P_1, \ldots, P_N$ with collections
\[
([v_{ij} : v_{ji}]) \in \prod_{1 \leq i < j \leq N} \bP(P_i \times P_j) = \E
\]
that satisfy $\| v_{ij} \|^2 = \| v_{ji} \|^2$ for all pairs. 
\end{remark}
By definition \emph{valid} common lines data is a collection $([v_{ij} : v_{ji}])$ of common lines data that has $N$ generic realizing frames $F_1, \ldots, F_N$. In coordinates this means that the associated embeddings bring together the common line pairs, i.e. for all $1 \leq i < j \leq N$, and for any representative $(v_{ij}, v_{ji})$, we have
\[
\iota_i(v_{ij}) = \iota_j(v_{ji}).
\]

\subsection{Necessary and sufficient conditions}
\label{sect-eqns}
In this section we derive equations and inequalities that are necessary and sufficient for common lines data $([v_{ij} : v_{ji}])$ to be valid. We first discuss necessary conditions. Recall from [Section~\ref{sect-results}] that for any triple of indices $i, j, k$ the angles between the common line pairs $[v_{ij} : v_{ji}]$, $[v_{ik} : v_{ki}]$ and $[v_{jk} : v_{kj}]$ determine a spherical triangle [Figure~\ref{fig-vanheel-tri}], and so these angles must satisfy the spherical triangle inequalities. The spherical triangle inequalities state that a non-degenerate spherical triangle of edge lengths $\alpha$, $\beta$ and $\gamma$, all in $(0, \pi)$, exists if and only if
\begin{equation}
\label{eqn-sphere}
\begin{aligned}
\beta + \gamma &> \alpha, \\
\alpha + \gamma &> \beta, \\
\alpha + \beta &> \gamma, \\
\alpha + \beta + \gamma &< 2\pi.
\end{aligned}
\end{equation}
\begin{definition}
Fix common lines data $([v_{ij} : v_{ji}]) \in \E$ and a triple of indices $(i,j,k)$. Choose representatives $(v_{ij}, v_{ji})$, $(v_{ik}, v_{ki})$ and $(v_{jk}, v_{kj})$. Then we write
\[
\alpha_{ijk} = \cos^{-1}\left( \frac{v_{ij} \cdot v_{ik}}{\| v_{ij} \| \| v_{ik} \|} \right), \quad \quad \beta_{ijk} = \cos^{-1}\left( \frac{v_{ji} \cdot v_{jk}}{\| v_{ji} \| \| v_{jk} \|} \right), \quad \quad \gamma_{ijk} = \cos^{-1}\left( \frac{v_{ki} \cdot v_{kj}}{\| v_{ki} \| \| v_{kj} \|} \right).
\]
\end{definition}
The angles $\alpha_{ijk}, \beta_{ijk}$ and $\gamma_{ijk}$ depend on the representatives we have chosen, however whether or not the spherical triangle inequalities [Equation~\ref{eqn-sphere}] are satisfied is independent of this choice. Thus the following definition makes sense:
\begin{definition}
Fix common lines data $([v_{ij} : v_{ji}]) \in \E$ and a triple of indices $(i,j,k)$. We say $(i,j,k)$ \uline{strictly satisfies the triangle inequalities} if, for any choice of representatives of $[v_{ij} : v_{ji}]$, $[v_{ik} : v_{ki}]$ and $[v_{jk} : v_{kj}]$, the angles $\alpha_{ijk}, \beta_{ijk}, \gamma_{ijk}$ satisfy [Equation~\ref{eqn-sphere}]. 
\end{definition}
This definition allows us to state our first result.
\begin{proposition}
\label{prop-three}
Suppose $([v_{ij} : v_{ji}]) \in \E$ is common lines data, and fix a triple $(i,j,k)$ of indices. Suppose that $(i,j,k)$ strictly satisfies the spherical triangle inequalities. Then there exist generic frames $F_i, F_j, F_k$ that realize the common line pairs $[v_{ij} : v_{ji}]$, $[v_{ik} : v_{ki}]$ and $[v_{jk} :v_{kj}]$. Moreover if $G_i, G_j, G_k$ are another set of frames that realize these same pairs, then there is a rotation $A \in \dO(3)$ such that $A$ maps $(F_i, F_j, F_k) \mapsto (G_i, G_j, G_k)$.
\end{proposition}
For a proof of this proposition, see [Appendix~\ref{proof-prop-three}]. 
\\

\begin{figure}
\begin{center}
\includegraphics[width=2in]{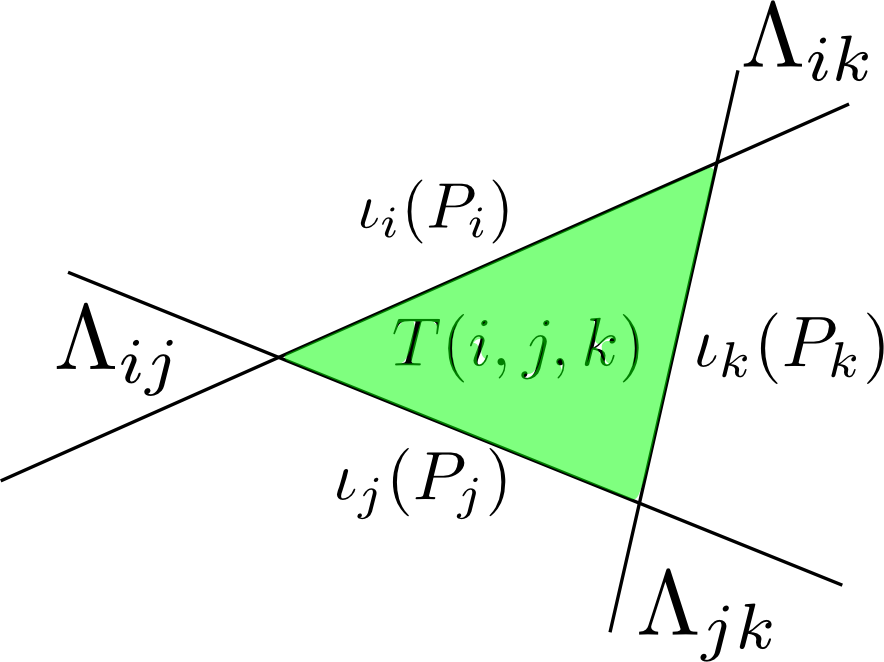}
\end{center}
\caption{Triangle $T(i,j,k)$ on the surface of the sphere.}
\label{fig-tri-ijk}
\end{figure}

This proposition is a necessary and sufficient condition for realizing frames to exist for a triple $(i,j,k)$, and so we have obtained necessary and sufficient conditions for $N = 3$. For $N > 3$, this proposition states that \emph{each} triple of indices $(i,j,k)$ must satisfy the spherical triangle inequality, but this condition is no longer sufficient.
\begin{example}
\label{ex-loc}
Consider the common lines data for $P_1, P_2, P_3, P_4$ given by
\begin{alignat*}{2}
(v_{12}, v_{13}, v_{14}) &= ( \left[ 1, \, 0 \right]^T, \left[ \sqrt{2}/2, \, \sqrt{2}/2 \right]^T, \left[ 0, \, 1 \right]^T ),& \quad (v_{21}, v_{23}, v_{24}) &= ( \left[ 1, \, 0 \right]^T, \left[ \sqrt{2}/2, \, \sqrt{2}/2 \right]^T, \left[ 0, \, 1 \right]^T ), \\ 
(v_{31}, v_{32}, v_{34}) &= ( \left[ 1, \, 0 \right]^T, \left[ \sqrt{2}/2, \, \sqrt{2}/2 \right]^T, \left[ 0, \, 1 \right]^T ),& \quad (v_{41}, v_{42}, v_{43}) &= ( \left[ 1, \, 0 \right]^T, \left[ \sqrt{2}/2, \, \sqrt{2}/2 \right]^T, \left[ 0, \, 1 \right]^T ).
\end{alignat*}
The angles between these common lines are given by
\begin{alignat*}{2}
(\alpha_{123}, \beta_{123}, \gamma_{123}) &= \left(\frac{\pi}{4},\frac{\pi}{4},\frac{\pi}{4}\right), \qquad & (\alpha_{124}, \beta_{124}, \gamma_{124}) &= \left(\frac{\pi}{2},\frac{\pi}{2},\frac{\pi}{4}\right), \\
(\alpha_{134}, \beta_{134}, \gamma_{134}) &= \left(\frac{\pi}{4},\frac{\pi}{2},\frac{\pi}{2}\right), \qquad & (\alpha_{234},\beta_{234},\gamma_{234}) &= \left(\frac{\pi}{4},\frac{\pi}{4},\frac{\pi}{4}\right).
\end{alignat*}
Observe that each of these triples satisfies the spherical triangle inequality. However, this data cannot be realized by frames $F_1, F_2, F_3$ and $F_4$. To see why, suppose such frames existed and, for each pair $i, j$, set $\Lambda_{ij} = \iota_i(v_{ij}) = \iota_j(v_{ji})$. The points $\Lambda_{12}, \Lambda_{13}, \Lambda_{23}$ determine a spherical triangle with edge lengths $(\alpha_{123}, \beta_{123}, \gamma_{123})$ [Figure~\ref{fig-loc-violate-1}], and the angle of this spherical triangle at the vertex between edges $\alpha_{123}$, $\beta_{123}$ is exactly the angle $\theta_{12}$ between the planes $\iota_1(P_1)$ and $\iota_2(P_2)$. From the spherical law of cosines, we can compute this angle:
\[
\cos \theta_{12} = \frac{\cos \gamma_{123} - \cos \alpha_{123} \cos \beta_{123}}{\sin \alpha_{123} \sin \beta_{123}} = \sqrt{2} - 1.
\]
Similarly, the points $\Lambda_{12}, \Lambda_{14}$ and $\Lambda_{24}$ determine a spherical triangle with edge lengths $(\alpha_{124}, \beta_{124}, \gamma_{124})$ [Figure~\ref{fig-loc-violate-2}], and the angle of this triangle between edges $\alpha_{124}$ and $\beta_{124}$ is again the angle $\theta_{12}$ between the planes $\iota_1(P_1)$ and $\iota_2(P_2)$. However, in this triangle we have $\cos \theta_{12} = \sqrt{2}/2$, which is a contradiction.
\begin{figure}[H]
\begin{center}
\begin{subfigure}[b]{0.45\textwidth}
\centering\includegraphics[width=2in]{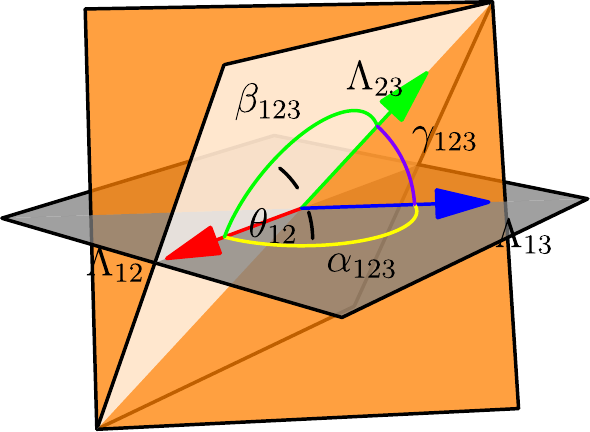}
\caption{Using $[v_{12}:v_{21}]$, $[v_{13} :v_{31}]$, $[v_{23}: v_{32}]$.}
\label{fig-loc-violate-1}
\end{subfigure}
\begin{subfigure}[b]{0.45\textwidth}
\centering\includegraphics[width=2in]{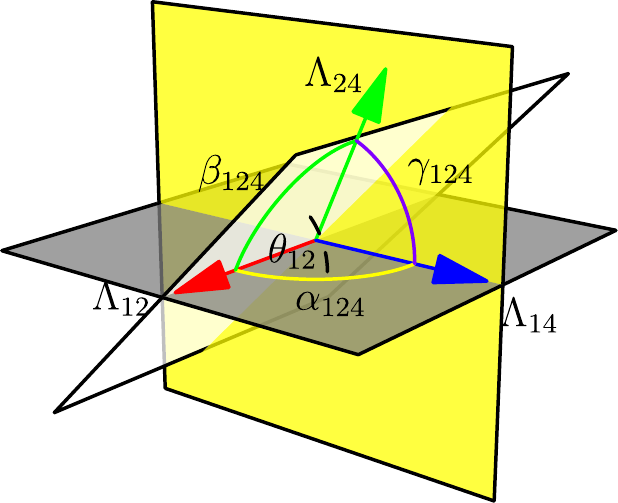}
\caption{Using $[v_{12}:v_{21}]$, $[v_{14} : v_{41}]$, $[v_{24} : v_{42}]$.}
\label{fig-loc-violate-2}
\end{subfigure}
\end{center}
\caption{Inconsistent reconstruction from invalid common lines data.}
\label{fig-loc-violate}
\end{figure}
\end{example}
We now discuss one explanation for why the contradiction in [Example~\ref{ex-loc}] arose that will lead us to a necessary and sufficient condition for reconstruction when $N > 3$. Suppose the frames $F_1, \ldots, F_N$ realize the common lines data $([v_{ij} : v_{ji}]) \in \E$, and choose unit vector representatives $(v_{ij}, v_{ji})$ for all the common line pairs. If we consider the intersection of the embedded planes $\iota_i(P_i)$ with the unit sphere in $\bR^3$, we obtain $N$ geodesic arcs. Each pair of these arcs has a distinguished point of intersection $\iota_i(v_{ij}) = \iota_j(v_{ji})$ which we denote by $\Lambda_{ij}$. Denote by $T(i,j,k)$ the triangle obtained by taking $\Lambda_{ij}$, $\Lambda_{ik}$ and $\Lambda_{jk}$ as vertices [Figure~\ref{fig-tri-ijk}]. 

\begin{figure}[h]
\begin{center}
\begin{subfigure}[b]{0.4\textwidth}
\centering\includegraphics[width=2.5in]{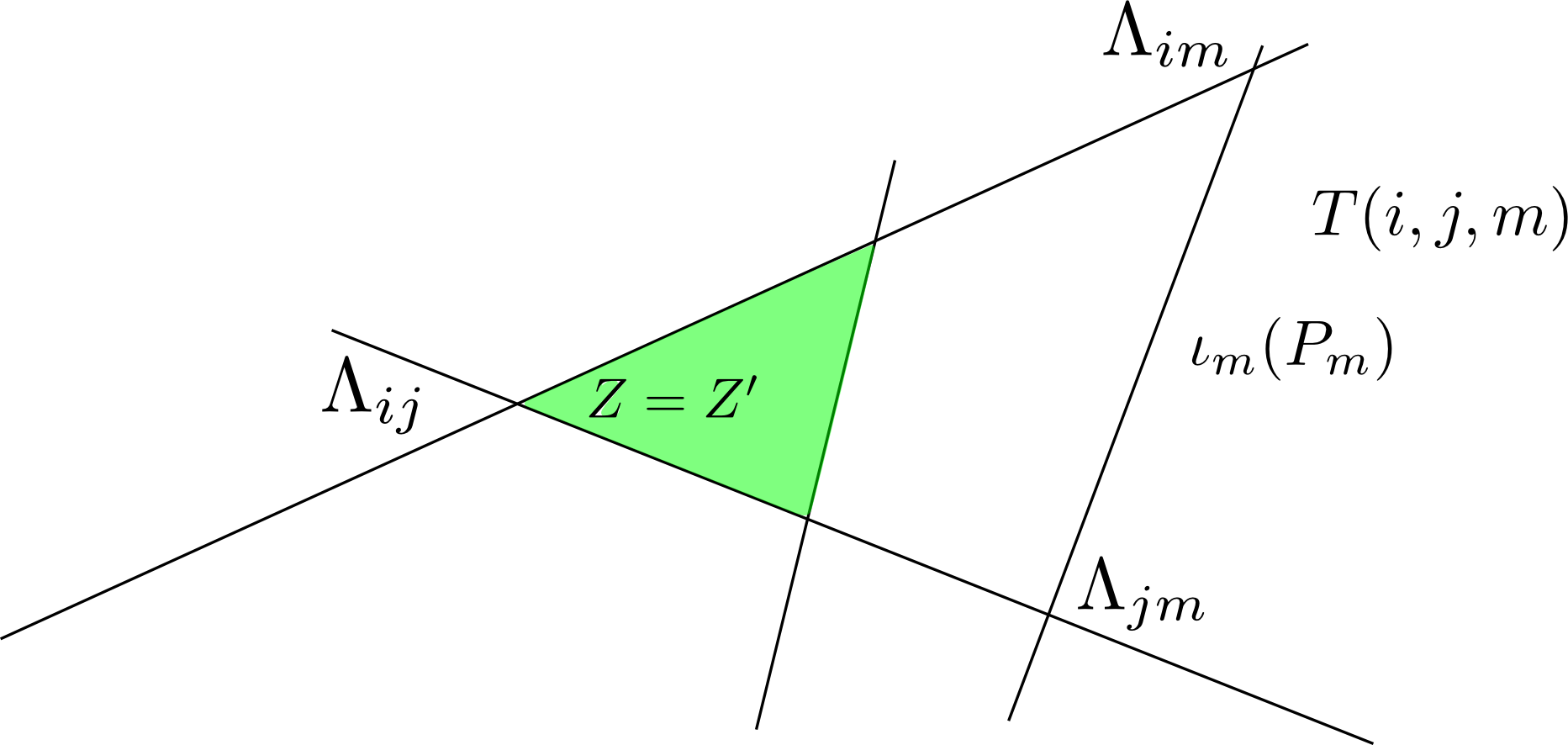}
\end{subfigure}
\qquad \qquad
\begin{subfigure}[b]{0.4\textwidth}
\centering\includegraphics[width=2.5in]{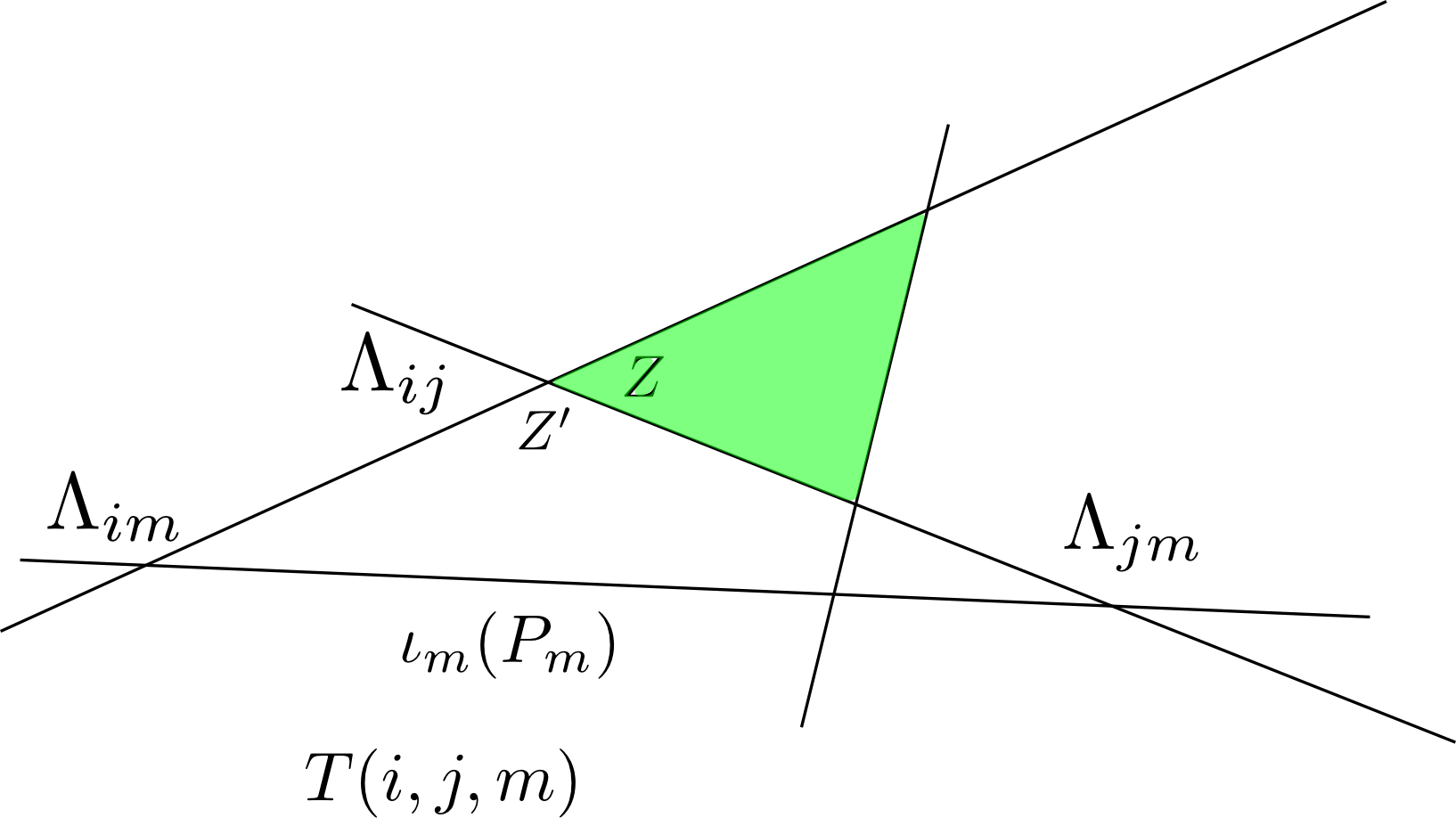}
\end{subfigure}
\end{center}
\caption{$T(i,j,m)$ shares edges with $T(i,j,k)$, so either $Z = Z'$ or $Z = \pi - Z'$.}
\label{fig-tri-ijm}
\end{figure}

Consider the second triangle $T(i,j,m)$ [Figure~\ref{fig-tri-ijm}]. The two triangles $T(i,j,k)$ and $T(i,j,m)$ share a vertex, $\Lambda_{ij}$, and the edges of \emph{both} triangles at this vertex lie in $\iota_i(P_i)$ and $\iota_j(P_j)$. It follows that the angle $Z$ in $T(i,j,k)$ and $Z'$ in $T(i,j,m)$ at this common vertex must be compatible: the angles are either the same, or supplementary, depending on the arrangement of the vertices. We can express this requirement in terms of the common lines data by using the spherical law of cosines:
\[
(\cos \gamma_{123} - \cos \alpha_{123} \cos \beta_{123}) \sin \alpha_{124} \sin \beta_{124} = \sigma(\cos \gamma_{124} - \cos \alpha_{124} \cos \beta_{124})\sin \alpha_{123} \sin \beta_{123},
\]
where $\sigma$ determines whether $Z = Z'$ or $Z = \pi - Z'$. In this light, the contradiction in [Example~\ref{ex-loc}] arose because the angles at $\Lambda_{12}$ in $T(1,2,3)$ and $T(1,2,4)$ were not compatible.

The spherical law of cosines compatibility described above is necessary for such a system of triangles $T(i,j,k)$ constructed from $F_1, \ldots, F_N$ to exist, and we will see it is sufficient as well. We first show that if this law of cosines compatibility between $(i,j,k)$ and $(i,j,m)$ is satisfied, then we can glue together reconstructions of these triples in a compatible fashion.
\newpage
\begin{lemma}
\label{lem-compat}
Suppose $([v_{ij} : v_{ji}]) \in \E$ is common lines data, and fix triples $(i,j,k)$ and $(i,j,m)$ that strictly satisfy the spherical triangle inequalities. Consider the spherical law of cosines compatibility
\begin{equation}
\label{eqn-loc}
\begin{aligned}
L_{ijk,ijm} = (v_{ki} \cdot v_{kj} - (v_{ij} \cdot v_{ik})&(v_{ji} \cdot v_{jk})) \mid \det \left[ v_{ij}, v_{im} \right] \det \left[ v_{ji}, v_{jm} \right] \mid - \\ & \sigma (v_{mi} \cdot v_{mj} - (v_{ij} \cdot v_{im})(v_{ji} \cdot v_{jm})) \mid \det \left[ v_{ij}, v_{ik} \right] \det \left[ v_{ji}, v_{jk} \right] \mid,
\end{aligned}
\end{equation}
where
\[
\sigma = \sign(\det \left[ v_{ij}, v_{ik} \right] \det \left[ v_{ij}, v_{im} \right ] \det \left[ v_{ji}, v_{jk} \right] \det \left[ v_{ji}, v_{jm} \right ]).
\]
Suppose that $L_{ijk,ijm} = 0$. Then, if $F_i, F_j, F_k$ are any realizing frames for $(i,j,k)$, and $G_i, G_j, G_m$ are any realizing frames for $(i,j,m)$, then there exists a unique rotation $A \in \dO(3)$ such that $AF_i = G_i$ and $AF_j = G_j$.
\end{lemma}
For a proof, see [Appendix~\ref{proof-lem-compat}]. 
\\

We now can show that the law of cosines compatibility is sufficient for reconstruction.
\begin{theorem}
\label{thm-reconstruct}
Suppose $([v_{ij} : v_{ji}]) \in \E$ strictly satisfies all spherical triangle inequalities and all spherical law of cosines compatibilities. Then there exist generic frames $F_1, \ldots, F_N$, unique up to isometry in $\dO(3)$, realizing $[v_{ij} : v_{ji}]$. 
\end{theorem}
For a proof, see [Appendix~\ref{proof-thm-reconstruct}].

\section{Geometry of valid common lines}
\label{sect-geom}
In this section we use the necessary and sufficient conditions derived in [Section~\ref{sect-eqns}] for valid common lines data to deduce some geometric properties about the set $C_N$ of all valid common lines. The main result in this section is that $C_N$ is homeomorphic to the space of generic frames, up to $\dO(3)$. In particular, this implies that the dimension of $C_N$ is $3N - 3$.
 
We first explicitly describe how to obtain valid common lines from a set of generic realizing frames $F_1, \ldots, F_N$ as in [Section~\ref{sect-common}]. For each pair $i,j$, choose a vector $\Lambda_{ij}$ in the one dimensional vector space $\iota_i(P_i) \cap \iota_j(P_j)$. Since $\bR^3$ has a canonical structure of an inner product space, we have the corresponding orthogonal projections $\iota^T_i: \bR^3 \to P_i$ and $\iota^T_j: \bR^3 \to P_j$. Consider the vectors
\[
(v_{ij}, v_{ji}) = (\iota^T_i(\Lambda_{ij}), \iota^T_j(\Lambda_{ij})) \in P_i \times P_j.
\]
By construction the pair $[v_{ij} : v_{ji}] = [ x_{ij} : y_{ij} : x_{ji} : y_{ji}]$ is a common line pair realized by the frames $F_i$ and $F_j$. In coordinates, we have
\begin{equation}
\label{eqn-common-coords}
x_{ij}a_i + y_{ij}b_i = \Lambda_{ij} = x_{ji}a_j + y_{ji}b_j.
\end{equation}
Repeating this process for all pairs $1 \leq i < j \leq N$, we obtain valid common lines data $([v_{ij} : v_{ji}]) \in C_N$ that is realized by $F_1, \ldots, F_N$. This algorithmically gives a map $\cG \to C_N$, where $\cG$ is the subset of $N$ generic frames in $\cF^N$. It will be useful to express this function via explicit polynomial mappings. We first describe a set of coordinates on the Grassmannian $\Gr(3,2N)$, whose points are the three dimensional subspaces of $\bR^{2N}$. 

\subsection{Grassmannian \& Pl\"ucker coordinates}
\label{sect-grass}
If $W \subset \bR^{2N}$ is a three dimensional subspace of $\bR^{2N}$, and we choose a basis $w_1, w_2, w_3 \in \bR^{2N}$ for $W$, we can represent the point in $\Gr(3,2N)$ corresponding to $W$ by the vector of all $3 \times 3$ minors of the $3 \times 2N$ matrix
\[
\left[ w_1, w_2, w_3 \right ]^T.
\]
These minors are the \emph{Pl\"ucker coordinates} of $W$. If we choose a different basis for $W$, the vector of $3 \times 3$ minors will only change by a non-zero scalar. Since Pl\"ucker coordinates are well defined up to scaling, the Grassmannian $\Gr(3,2N)$ is a subvariety of the projective space $\bP^{\binom{2N}{3} - 1}$.

Given a collection of $N$ frames $F_1, \ldots, F_N$, we can form the $3 \times 2N$ matrix
\[
F_\bullet = [F_1\, \ldots \, F_N] = [a_1, b_1, \ldots, a_N, b_N].
\]
We consider the rational map $\rho: \cF^N \dashrightarrow \Gr(3,2N)$ that takes a collection of frames $F_1, \ldots, F_N$ to the Pl\"ucker coordinates of $F_\bullet$. A \emph{rational map} is a map that is defined almost everywhere in the domain. In this case, $\rho$ is not defined if the rank of $F_\bullet$ is $\leq 2$, since, in this case, the rows of $F_\bullet$ do not determine a three dimensional subspace of $\bR^{2N}$. 
\subsection{Pl\"ucker coordinates for common lines}
As described above, given a pair of frames $F_i, F_j$ we can compute the associated common line pair $[v_{ij} : v_{ji}] \in C_N$ by choosing \emph{any} vector $\Lambda_{ij}$ in $\iota_i(P_i) \cap \iota_j(P_j)$. In particular, we can choose $\Lambda_{ij} = (a_i \times b_i) \times (a_j \times b_j)$, where $\times$ is the standard vector cross product on $\bR^3$. Then, the following identity from vector algebra, called the vector quadruple product, expresses $\Lambda_{ij}$ in terms of the frames $F_i$ and $F_j$:
\[
\det \left[ a_j, b_j, a_i \right] b_i - \det \left[ a_j, b_j, b_i \right] a_i = (a_i \times b_i) \times (a_j \times b_j) = \det \left [ a_i, b_i, b_j \right] a_j - \det \left[ a_i, b_i, a_j \right] b_j.
\]
Comparing this with [Equation~\ref{eqn-common-coords}], we see that the coordinates of the common line pair $[v_{ij} : v_{ji}]$ are given by determinants of certain $3 \times 3$ matrices. Explicitly, we have
\[
v_{ij} = \left [ \begin{array}{r} -\det \left[ a_j, b_j, b_i \right ]
    \\ \det \left[ a_j, b_j, a_i \right] \end{array} \right], \quad
\quad v_{ji} = \left[ \begin{array}{r} \det \left[ a_i, b_i, b_j
    \right ] \\ -\det \left[ a_i, b_i, a_j \right] \end{array} \right].
\]
Observe that these $3 \times 3$ determinants are certain $3 \times 3$ minors of the matrix $F_\bullet$. The minors that appear are those that belong to only two frames $F_i$ and $F_j$: in other words, any three of $\{ a_i, b_i, a_j, b_j \}$. The minors \emph{not} appearing as coordinates of a common line pair are those that choose three columns from three distinct frames:
\begin{equation}
\label{eqn-bad-minor}
\det \left [ \{ a_i, b_i \}, \{a_j, b_j \}, \{ a_k, b_k \} \right ].
\end{equation}
Thus, the coordinates on the Grassmannian $\Gr(3,2N)$ are the common line coordinates, together with these ``bad'' minors [Equation~\ref{eqn-bad-minor}]. If we consider the projection where we discard the ``bad'' minors, we obtain the map
\[
\Gr(3,2N) \dashrightarrow \prod_{1 \leq i < j \leq N} \bP(P_i \times P_j) = \E.
\]
Explicitly, for each $i, j$, this projection maps
\[
\left [ \ldots : -\det \left[ a_j, b_j, b_i \right ] : \det \left[ a_j, b_j, a_i \right] : \det \left[ a_i, b_i, b_j \right ] : -\det \left[ a_i, b_i, a_j \right] : \ldots \right ] \mapsto [v_{ij} : v_{ji}].
\]
Note that this rational map is not defined whenever the four $3 \times 3$ minors appearing in the common line pair $[v_{ij} : v_{ji}]$ simultaneously vanish. This cannot happen with generic frames, so this projection is an honest morphism when restricted to $\rho(\cG) \subset \Gr(3,2N)$. The image of this map is the set of valid common lines, and is in fact a homeomorphism.
\begin{theorem}
\label{thm-homeo}
The restriction $\pi$ of the projection $\Gr(3,2N) \dashrightarrow \E$ to $\rho(\cG) \subset \Gr(3,2N)$ is a homeomorphism onto $C_N$.
\end{theorem}
For a proof, see [Appendix~\ref{proof-thm-homeo}]. 
\\

This theorem means that we can identify $C_N$ with $\rho(\cG)$, which is an open subset of the Grassmannian. As we discussed above, the point $\rho(F_\bullet) \in \Gr(3,2N)$ only determines the row space of the matrix $F_\bullet = [F_1, \, \ldots, \, F_N]$. A different basis for this row space is given by multiplying $F_\bullet$ on the left by a matrix $A$ in $\dO(3)$, or, equivalently, by the following action
\[
A \cdot (F_1, \ldots, F_N) = (AF_1, \ldots, AF_N).
\] 
This is the \emph{diagonal action} of $\dO(3)$ on the space of frames $\cF^N$. We observe that this rotational ambiguity is the only difference between the space of frames and the Pl\"ucker embedding of these frames.
\begin{theorem}
\label{thm-homeo-frames}
The collection of generic frames $\cG \subset \cF^N$ is homeomorphic to $\rho(\cG) \times \dO(3)$.
\end{theorem}
For a proof, see [Appendix~\ref{proof-thm-homeo-frames}]. 
\\

Thus, we obtain the remainder of our results.
\begin{corollary}
The set $C_N$ of valid common lines is homeomorphic to the quotient of the set of generic frames $\cG$ by the group of rotations $\dO(3)$ acting diagonally.
\end{corollary}
In other words, we have recovered the fact that common lines data only determines its realizing frames up to $\dO(3)$: it is because we can identify common lines data with elements of $\rho(\cG) \subset \Gr(3,2N)$, and points in this space determine frames up to global rotation.
\begin{corollary}
\label{cor-dim}
The dimension of $C_N$ as a semi-algebraic set is $3N - 3$.
\end{corollary}
For a proof, see [Appendix~\ref{proof-cor-dim}].

\subsection{Defining Polynomials}
In [Section~\ref{sect-eqns}] we derived the defining equations for $C_N$ in terms of spherical geometry. For the benefit of the reader we now explicitly describe these conditions as multi-homogeneous polynomials in the variables $([v_{ij} : v_{ji}])$.

Suppose $([v_{ij} : v_{ji}]) \in \E$ is fixed, and that $\| v_{ij} \|^2 = \| v_{ji} \|^2$ for all $1 \leq i < j \leq N$. The spherical triangle inequalities for the common line pairs $[v_{ij} : v_{ji}]$, $[v_{ik} : v_{ki}]$ and $[v_{jk} : v_{kj}]$, [Equation~\ref{eqn-sphere}], are equivalent, see~\cite[Lemma 2.2]{buckman}, to
\[
\| v_{ij} \|^2\|v_{ik}\|^2\|v_{jk}\|^2 - \| v_{jk} \|^2 (v_{ij} \cdot v_{ik})^2 - \| v_{ik} \|^2(v_{ji} \cdot v_{jk})^2 - \| v_{ij} \|^2(v_{ki} \cdot v_{kj})^2 + 2(v_{ij} \cdot v_{ik})(v_{ji} \cdot v_{jk})(v_{ki} \cdot v_{kj}) > 0.
\]
To express the spherical law of cosines compatibilities $L_{ijk,ijm}$ [Equation~\ref{eqn-loc}], set 
\begin{align*}
a &= (\| v_{ij} \|^2 (v_{ki} \cdot v_{kj}) - (v_{ij} \cdot v_{ik})(v_{ji} \cdot v_{jk})), \\
b &= ( \| v_{ij} \|^2 (v_{mi} \cdot v_{mj}) - (v_{ij} \cdot v_{im})(v_{ji} \cdot v_{jm})), \\
d_1 &= \det[ v_{ij}, v_{im} ]\det[ v_{ji}, v_{jm} ], \\
d_2 &= \det[v_{ij}, v_{ik} ]\det [ v_{ji}, v_{jk} ].
\end{align*}
Then $L_{ijk,ijm} = 0$ if and only if
\[
a^2d_1^2 - 2d_1d_2ab + b^2d_2^2 = 0.
\]

Thus, the set $C_N$ is defined as a semi-algebraic subset of $\E$ by the following equations and inequalities:
\begin{enumerate}
\item The $\binom{N}{2}$ equations $\| v_{ij} \|^2 = \| v_{ji} \|^2$, see [Section~\ref{sect-coords}].
\item For each of the $\binom{N}{3}$ triples $(i,j,k)$ the spherical triangle inequality, see [Proposition~\ref{prop-three}].
\item For each of the $6\binom{N}{4}$ ways to choose two triples of distinct indices $(i,j,k), (i,j,m)$ the spherical law of cosines compatibility, see [Lemma~\ref{lem-compat}].
\end{enumerate}

\section{Acknowledgements}
The author is greatly thankful to Shamgar Gurevich, for initially suggesting the cryo-EM problem and for his continued support, as well as to Bernd Sturmfels, who suggested studying defining equations in cryo-EM, provided helpful guidance, and invited the author to the Mathematical Sciences Research Institute (MSRI) in Berkeley, California. Much of this work took place at MSRI during the spring of 2013, and the author greatly appreciates helpful technical discussions with Luke Oeding, Kristian Ranestad, Yoel Shkolnisky, Amit Singer and Frank Sottile. The author's visit to MSRI was supported by the National Science Foundation under grants 0838210 and 0932078.

\appendix
\section{Proofs}
\label{sect-proofs}
\begin{proof}[\bf Proof of Proposition~\ref{prop-three}]
\label{proof-prop-three}
Fix representatives $(v_{ij}, v_{ji}), (v_{ik}, v_{ki})$ and $(v_{jk}, v_{kj})$. Since the lengths $\alpha_{ijk}$, $\beta_{ijk}$, $\gamma_{ijk}$ strictly satisfy the triangle inequalities, there is a non-degenerate spherical triangle with these edge lengths. Denote the vertex of this triangle opposite the edge of length $\alpha_{ijk}$ by $\Lambda_{jk}$, the vertex opposite the edge $\beta_{ijk}$ by $\Lambda_{ik}$ and the vertex opposite the edge $\gamma_{ijk}$ by $\Lambda_{ij}$. Since this triangle is non-degenerate, we know that $\Lambda_{ij}, \Lambda_{ik}$ and $\Lambda_{jk}$ are linearly independent. Thus we have embeddings $\iota_i, \iota_j, \iota_k$, given by
\begin{alignat*}{3}
\iota_i: P_i &\hookrightarrow \bR^3, \quad & \iota_j: P_j &\hookrightarrow \bR^3, \quad & \iota_k: P_k &\hookrightarrow \bR^3, \\
v_{ij} &\mapsto \Lambda_{ij}, \quad & v_{ji} &\mapsto \Lambda_{ij}, \quad & v_{ki} &\mapsto \Lambda_{ik}, \\
v_{ik} &\mapsto \Lambda_{ik}, \quad & v_{jk} &\mapsto \Lambda_{jk}, \quad & v_{kj} &\mapsto \Lambda_{jk}.
\end{alignat*}
Observe that these embeddings are isometric by construction, and so $F_i = (\iota_i(x), \iota_i(y)), F_j = (\iota_j(x),\iota_j(y))$ and $F_k = (\iota_k(x),\iota_k(y))$ are frames. Since $\Lambda_{ij}, \Lambda_{ik}, \Lambda_{jk}$ are vertices of a non-degenerate spherical triangle, these three frames are in generic position. Moreover, by construction we have
\[
\iota_i(v_{ij}) = \iota_j(v_{ji}), \qquad \iota_i(v_{ik}) = \iota_k(v_{ki}), \qquad \iota_j(v_{jk}) = \iota_k(v_{kj}),
\]
and so $F_i, F_j$ and $F_k$ realize the required common line pairs.

Now, suppose $G_i, G_j, G_k$ realize the common line pairs $[v_{ij} : v_{ji}], [v_{ik} : v_{ki}]$ and $[v_{jk} : v_{kj}]$. Let $\iota^G_i$, $\iota^G_j$ and $\iota^G_k$ be the embeddings corresponding to these frames, and set $\Lambda^G_{ij} = \iota^G_i(v_{ij})$, $\Lambda^G_{ik} = \iota^G_i(v_{ik})$ and $\Lambda^G_{jk} = \iota^G_j(v_{jk})$. Since $(i,j,k)$ strictly satisfies the triangle inequalities, these three vectors are linearly independent and thus define a unique spherical triangle with edge lengths $(\alpha_{ijk}, \beta_{ijk}, \gamma_{ijk})$. This triangle is congruent to the triangle with vertices $\Lambda_{ij}$, $\Lambda_{ik}$, $\Lambda_{jk}$ constructed above, and so there exists an isometry $A \in \dO(3)$ that maps $\Lambda^G_{ij} \mapsto \Lambda_{ij}, \Lambda^G_{ik} \mapsto \Lambda_{ik}$ and $\Lambda^G_{jk} \mapsto \Lambda_{jk}$, and thus maps $(G_i, G_j, G_k) \mapsto (F_i, F_j, F_k)$.
\end{proof}

\begin{proof}[\bf Proof of Lemma~\ref{lem-compat}]
\label{proof-lem-compat}
Fix representatives $(v_{ij}, v_{ji}), (v_{ik}, v_{ki}), (v_{jk}, v_{kj}), (v_{im}, v_{mi})$ and $(v_{jm}, v_{mj})$. Let $\iota^F_i$, $\iota^F_j$, and $\iota^F_k$ be the embeddings corresponding to $F_i$, $F_j$, $F_k$, and let $\iota^G_i$, $\iota^G_j$, $\iota^G_m$ be the embeddings corresponding to $G_i$, $G_j$, $G_m$. Consider the embedded common lines $\Lambda^F_{ij} = \iota^F_i(v_{ij})$, $\Lambda^F_{ik} = \iota^F_i(v_{ik})$ and $\Lambda^F_{jk} = \iota^F_j(v_{jk})$, and similarly let $\Lambda^G_{ij} = \iota^G_i(v_{ij})$, $\Lambda^G_{ik} = \iota^G_i(v_{ik})$ and $\Lambda^G_{jk} = \iota^G_j(v_{jk})$. Since the triple $(i,j,k)$ strictly satisfies the triangle inequality, we know that both these sets of vectors are linearly independent. Let $A: \bR^3 \to \bR^3$ be the map defined by
\[
\Lambda^F_{ij} \mapsto \Lambda^G_{ij}, \qquad \Lambda^F_{ik} \mapsto \Lambda^G_{ik}, \qquad \Lambda^F_{jk} \mapsto \Lambda^G_{jk}.
\]
We wish to show that $A$ is an isometry. Observe that 
\begin{align*}
\Lambda^F_{ij} \cdot \Lambda^F_{ik} = \iota^F_i(v_{ij}) \cdot \iota^F_i(v_{ik}) = v_{ij} &\cdot v_{ik} = \iota^G_i(v_{ij}) \cdot \iota^G_i(v_{ik}) = \Lambda^G_{ij} \cdot \Lambda^G_{ik}, \\
\Lambda^F_{ij} \cdot \Lambda^F_{jk} = \iota^F_j(v_{ji}) \cdot \iota^F_j(v_{jk}) = v_{ji} &\cdot v_{jk} = \iota^G_j(v_{ji}) \cdot \iota^G_j(v_{jk}) = \Lambda^G_{ij} \cdot \Lambda^G_{jk},
\end{align*}
so we only need to show that $\Lambda^F_{ik} \cdot \Lambda^F_{jk} = \Lambda^G_{ik} \cdot \Lambda^G_{jk}$.

Consider first the three points $\iota^G_i(v_{ij})$, $\iota^G_i(v_{im})$, $\iota^G_j(v_{jm})$ on the unit sphere in $\bR^3$. Since the triple $(i,j,m)$ strictly satisfies the triangle inequalities, these three points are the vertices of a unique spherical triangle. Let $Z$ be the angle of this triangle at the vertex $\Lambda^G_{ij}$.

Now, observe that the vectors $\Lambda^G_{ij}$ and $-\Lambda^G_{ij}$ cut the unit circle in the $\iota^G_i(P_i)$ plane into two semi-circles, and the points $\Lambda^G_{ik}$ and $\Lambda^G_{im}$ lie in these semi-circles. Furthermore, if $\sign(\det [ v_{ij}, v_{ik} ]) = \sign(\det [v_{ij}, v_{im}])$, then both of these vectors are in the same semi-circle. Similarly $\Lambda^G_{jk}$ and $\Lambda^G_{jm}$ lie in the same semi-circle in the $\iota^G_j(P_j)$ plane if $\sign(\det [ v_{ji}, v_{jk} ]) = \sign(\det [v_{ji}, v_{jm}])$. 

The points $\Lambda^G_{ij}, \Lambda^G_{ik}, \Lambda^G_{jk}$ define a unique spherical triangle. Let $Z'$ be the vertex of this triangle at $\Lambda^G_{ij}$. We claim that $Z'$ is either $Z$ or $\pi - Z$, depending on which semi-circles the points $\Lambda^G_{ik}$, $\Lambda^G_{jk}$, $\Lambda^G_{im}$, and $\Lambda^G_{jm}$ lie in. We have $Z' = Z$ if $\Lambda^G_{ik}, \Lambda^G_{im}$ are in the same semi-circle in $\iota^G_i(P_i)$ and $\Lambda^G_{jk},\Lambda^G_{jm}$ are in the same semi-circle in $\iota^G_j(P_j)$, or if $\Lambda^G_{ik}, \Lambda^G_{im}$ are in opposite semi-circles \emph{and} $\Lambda^G_{jk},\Lambda^G_{jm}$ are in opposite semi-circles. If one pair is in the same semi-circle, and the other pair is in opposite semi-circles, then we have $Z' = \pi - Z$. Equivalently, we have $Z = Z'$ if $\sigma = \sign(\det [ v_{ij}, v_{ik} ]\det [v_{ij}, v_{im}]\det [ v_{ji}, v_{jk} ]\det [v_{ji}, v_{jm}]) > 0$ and $Z' = \pi - Z$ if this product is negative.

Now, based on our discussion above, the cosine of the angle at vertex $\Lambda^G_{ij}$ in the triangle $\Lambda^G_{ij}$, $\Lambda^G_{ik}$, $\Lambda^G_{jk}$ is $\sigma \cos Z'$, so from the spherical law of cosines we obtain
\[
\Lambda^G_{ik} \cdot \Lambda^G_{jk} = (v_{ij} \cdot v_{ik})(v_{ji} \cdot v_{jk}) + \mid \det [v_{ij}, v_{ik}] \det [v_{ji}, v_{jk}] \mid \sigma \cos Z'.
\]
On the other hand, $L_{ijk,ijm} = 0$ implies that $\sigma \cos Z'$ is the angle at $\Lambda^F_{ij}$ in triangle $\Lambda^F_{ij}$, $\Lambda^F_{ik}$, $\Lambda^F_{jk}$, so again applying the law of cosines we have
\begin{align*}
\frac{\Lambda^G_{ik} \cdot \Lambda^G_{jk} - (v_{ij} \cdot v_{ik})(v_{ji} \cdot v_{jk})}{\mid \det [v_{ij}, v_{ik}] \det [v_{ji}, v_{jk}] \mid}&= \frac{v_{ik} \cdot v_{jk} - (v_{ij} \cdot v_{ik})(v_{ji} \cdot v_{jk})}{\mid \det [v_{ij}, v_{ik}] \det [v_{ji}, v_{jk}] \mid},
\end{align*}
and thus $\Lambda^G_{ik} \cdot \Lambda^G_{jk} = v_{ki} \cdot v_{ki} = \Lambda^F_{ik} \cdot \Lambda^F_{jk}$. We conclude that $A$ is an isometry, as desired.
\end{proof}

\begin{proof}[\bf Proof of Theorem~\ref{thm-reconstruct}]
\label{proof-thm-reconstruct}
Choose representatives $(v_{ij}, v_{ji})$ for all common line pairs, and suppose that $F_i$, $F_j$, $F_k$ and $G_i$, $G_j$, $G_m$ are two sets of realizing frames for their associated common lines. From [Lemma~\ref{lem-compat}] we know that there exists a unique rotation $A$ that maps $F_i \mapsto G_i$ and $F_j \mapsto G_j$. Observe that if $\det A = -1$ we can replace $G_i$, $G_j$, $G_m$ by $LG_i$, $LG_j$, $LG_m$, where $L$ is an arbitrary rotation with determinant $-1$. Note that $LG_i$, $LG_j$, $LG_m$ is still a set of realizing frames, but the compatibility morphism is now $L \circ A$ which has determinant $1$.

We proceed with reconstructing a set of realizing frames. By [Proposition~\ref{prop-three}], we first obtain frames $F_1$, $F_2$, $F_3$ from the triple $(1,2,3)$. For all remaining $i$, we choose a reconstruction $G_1$, $G_2$, $G_i$ from the triple $(1,2,i)$. By [Lemma~\ref{lem-compat}], there exists a unique map $A_i$ that maps $F_1 \mapsto G_1$ and $F_2 \mapsto G_2$, and, by our discussion above, we can assume that $A_i$ has determinant $1$. We set $F_i = A^{-1}_iG_i$.

We now need to check that these frames are realizing frames. We will write $\iota^F_i$, $\iota^F_j$ and $\iota^F_k$ for the embeddings determined by $F_i$, $F_j$, $F_k$, and similarly for the other sets of reconstructed frames. Thus we need to verify that $\iota^F_i(v_{ij}) = \iota^F_j(v_{ji})$. To this end, suppose that $F_i = A^{-1}_iG_i$ was reconstructed from $G_1$, $G_2$, $G_i$ and $F_j = A^{-1}_jD_j$ was reconstructed from $D_1$, $D_2$, $D_j$. The triple $(1,i,j)$ also strictly satisfies the triangle inequality, so we have generic realizing frames $H_1$, $H_i$, $H_j$. By [Lemma~\ref{lem-compat}] we obtain unique rotations $B_i: (G_1, G_i) \mapsto (H_1, H_i)$ and $B_j: (D_1, D_j) \mapsto (H_1, H_j)$.

We claim that $\det B_i = \det B_j$. This follows, since
\begin{align*}
\det \left[ \iota^G_1(v_{12}), \iota^G_1(v_{1i}), \iota^G_i(v_{i2}) \right ]&= (\det B_i) \left[ \iota^H_1(v_{12}), \iota^H_i(v_{i1}), \iota^H_i(v_{i2}) \right], \\
\det \left[ \iota^D_1(v_{12}), \iota^D_1(v_{1j}), \iota^D_j(v_{j2}) \right ]&= (\det B_j) \left[ \iota^H_1(v_{12}), \iota^H_j(v_{j1}), \iota^H_j(v_{j2}) \right].
\end{align*}
On the other hand, if $\sigma = \sign(\det[ v_{12}, v_{1i} ] \det[ v_{12}, v_{1j} ] \det[v_{21}, v_{2i}] \det[ v_{21}, v_{2j}])$, then
\begin{align*}
\det \left[ \iota^G_1(v_{12}), \iota^G_1(v_{1i}), \iota^G_i(v_{i2}) \right ] &= \sigma \det \left[ \iota^D_1(v_{12}), \iota^D_1(v_{1j}), \iota^D_j(v_{j2}) \right ], \\
\det \left[ \iota^H_1(v_{12}), \iota^H_i(v_{i1}), \iota^H_i(v_{i2}) \right] &= \sigma \det \left[ \iota^H_1(v_{12}), \iota^H_j(v_{j1}), \iota^H_j(v_{j2}) \right],
\end{align*}
and so $\det B_i = \det B_j$. Now, consider the diagram
\[
\xymatrix{
& (G_1, G_2, G_i) \ar[rd]^{B_i} & \\
(F_1, F_2, F_3) \ar[ru]^{A_i} \ar[rd]_{A_j} & & (H_1, H_i, H_j) \\
& (D_1, D_2, D_j) \ar[ru]_{B_j} & 
}
\]
Note that this diagram commutes, since both the top path and bottom path are morphisms in $\dO(3)$ of the same determinant that send $F_1 \mapsto H_1$. Then, since $H_1$, $H_i$, $H_j$ realize the common line pair $(v_{ij}, v_{ji})$, we have
\begin{align*}
\iota^F_i(v_{ij}) = A^{-1}_i \iota^G_i(v_{ij}) &= (A^{-1}_i \circ B^{-1}_i) \iota^H_i(v_{ij}) \\
&= (A^{-1}_i \circ B^{-1}_i) \iota^H_j(v_{ji}) \\
&= (A^{-1}_j \circ B^{-1}_j) \iota^H_j(v_{ji}) = A^{-1}_j \iota^G_j(v_{ji}) = \iota^F_j(v_{ji})
\end{align*} 
and thus $F_1, \ldots, F_N$ realize the common lines data $[v_{ij} : v_{ji}]$.

Finally, suppose that $F'_1, \ldots, F'_N$ is another collection of frames that realizes $[v_{ij} : v_{ji}]$. Choose three indices $(i,j,k)$, and consider the vectors $\Lambda'_{ij} = \iota^{F'}_i(v_{ij})$, $\Lambda'_{ik} = \iota^{F'}_i(v_{ik})$ and $\Lambda'_{jk} = \iota^{F'}_j(v_{jk})$. The angles between these vectors are given by $(\alpha_{ijk}, \beta_{ijk}, \gamma_{ijk})$. Similarly the angles between $\Lambda_{ij} = \iota^F_i(v_{ij})$, $\Lambda_{ik} = \iota^F_i(v_{ik})$ and $\Lambda_{jk} = \iota^F_j(v_{jk})$ are also given by $(\alpha_{ijk}, \beta_{ijk}, \gamma_{ijk})$, so there is an isometry $R_{ijk}$ that sends $(F'_i, F'_j, F'_k) \mapsto (F_i, F_j, F_k)$. Note that, $R_{ijk} = R_{ijm}$ for any $i,j,k,m$, since these are two isometries that agree on $F'_i$ and $F'_j$. This implies that $R_{ijk}(F'_m) = F_m$ for all $m$, and thus there is a single isometry $R: (F'_1, \ldots, F'_N) \mapsto (F_1, \ldots, F_N)$.
\end{proof}

\begin{proof}[\bf Proof of Theorem~\ref{thm-homeo}]
\label{proof-thm-homeo}
First, observe that the minors corresponding to a common line pair $[v_{ij} : v_{ji}]$ are non-zero for points in $\rho(\cG)$, since otherwise $F_i$ and $F_j$ would define the same plane. It follows that the restriction of the rational projection $\Gr(3,2N) \dashrightarrow \E$ to $\rho(\cG)$ yields an actual morphism. 

By our discussion above we know that $\pi(\rho(\cG)) = C_N$, so $\pi$ is onto $C_N$. It suffices to check that this projection is injective. This follows from [Theorem~\ref{thm-reconstruct}]. If $\pi(\rho(F_\bullet)) = \pi(\rho(G_\bullet))$, then we know that the realizing frames $F_\bullet$ and $G_\bullet$ are related by an isometry in $\dO(3)$. But then the rows of the matrices $F_\bullet$ and $G_\bullet$ define the same linear subspace, and so $\rho(F_\bullet) = \rho(G_\bullet)$.
\end{proof}

\begin{proof}[\bf Proof of Theorem~\ref{thm-homeo-frames}]
\label{proof-thm-homeo-frames}
For any $F_\bullet \in \cG$, for all indices $1 \leq i < j \leq N$, define $\Lambda^F_{ij} = (a_i \times b_i) \times (a_j \times b_j)$. Then, consider the map $\eta: \cG \to \rho(\cG) \times \dO(3)$, given by
\[
[a_1, b_1, \, \ldots \, a_N, b_N] \mapsto (\rho(F_\bullet), [a_1, b_1, \xi a_1 \times b_1]),
\]
where $\xi = \sign \det[ \Lambda^F_{12}, \Lambda^F_{13}, \Lambda^F_{23} ]$. This map is surjective so we just need to verify injectivity. Suppose that
\[
\eta(F_\bullet) = (\rho(F_\bullet), [a_1, b_1, \xi a_1 \times b_1]) = (\rho(F'_\bullet), [a'_1, b'_1, \xi' a'_1 \times b'_1]) = \eta(F'_\bullet).
\]
Since $\rho(F_\bullet) = \rho(F'_\bullet)$, we know that $R F_\bullet = F'_\bullet$, for some $R \in \dO(3)$. Furthermore, since $a_1 = a'_1$ and $b_1 = b'_1$, if we can show $\det(R) = 1$, then $R$ will be the identity.

To this end, observe that since $R F_i = F'_i$ for all $i$, we have,
\[
\Lambda^{F'}_{ij} = (R a_i \times R b_i) \times (R a_j \times
R b_j) = R \left( (a'_i \times b'_i) \times (a'_j \times b'_j)
\right ) = R \Lambda^F_{ij} \\
\]
for all $1 \leq i < j \leq N$. Then we have
\[
\det [ R\Lambda^F_{12}, R\Lambda^F_{13}, R\Lambda^F_{23} ] = \det R \det [ \Lambda^{F}_{12}, \Lambda^{F}_{13}, \Lambda^{F}_{23}] = \det[ \Lambda^{F'}_{12}, \Lambda^{F'}_{13}, \Lambda^{F'}_{23}],
\]
and, since $\xi = \xi'$, we must have $\det R = 1$. We conclude that $R$ is the identity, and thus $F_\bullet = F'_\bullet$.
\end{proof}

\begin{proof}[\bf Proof of Corollary~\ref{cor-dim}]
\label{proof-cor-dim}
The dimension of $\cG \subset \cF^N$ is the same as that of $\cF^N$ since $\cG$ is the complement of an algebraic hypersurface. The set of all frames $\cF^N$ is isomorphic to $\SO(3)^N$ and so has dimension $3N$. It follows that $\rho(\cG) \times \dO(3)$ also has dimension $3N$, and so $\rho(\cG) \cong C_N$ has dimension $3N - 3$.
\end{proof}

\bibliographystyle{alpha}
\bibliography{common_lines_geometry_open}

\end{document}